\documentclass[12pt,amsfonts]{article}
\usepackage{etex}
\usepackage{dsfont}
\usepackage{graphicx}
\usepackage{latexsym}
\usepackage{amssymb}
\usepackage{amsmath}
\usepackage{enumerate}
\usepackage{stmaryrd}
\usepackage{mathrsfs}
\usepackage{layout}
\usepackage{cases}
\usepackage{eufrak}
\usepackage{euscript}
\usepackage{pstricks, pst-plot, pst-node, pst-tree}
\usepackage{auto-pst-pdf}
\usepackage{bm}   

\usepackage{epsfig}
\usepackage{pst-grad} 
\usepackage{pst-plot} 

\usepackage{amsfonts}
\usepackage{float}

\usepackage[all]{xy}

\newtheorem{prop}{Proposition}[section]
\newtheorem{lemma}{Lemma}[section]
\newtheorem{definition}{Definition}[section]
\newtheorem{corollary}{Corollary}[section]
\newtheorem{theorem}{Theorem}[section]
\newtheorem{remark}{Remark}[section]

\newtheorem{assumption}{Assumption}[section]

\numberwithin{equation}{section}

\newcounter{alphasect}
\def\alphainsection{0}

\let\oldsection=\section
\def\section{%
  \ifnum\alphainsection=1%
    \addtocounter{alphasect}{1}
  \fi%
\oldsection}%

\renewcommand\thesection{%
  \ifnum\alphainsection=1%
    \Alph{alphasect}%
  \else%
    \arabic{section}%
  \fi%
}%

\newcommand{\E}{\mathbb{E}}

\def\real{{\mathord{{\rm I\kern-2.8pt R}}}}        
\def\inte{{\mathord{{\rm I\kern-2.8pt N}}}}

\def\sZZ{{\rm Z\kern-2.8ptem{}Z}}

\def\z{{\mathchoice
  {\sZZ}
  {\sZZ}
  {\rm Z\kern-0.30em{}Z}
  {\rm Z\kern-0.25em{}Z} }}
\def\sQQ{{\kern 0.27em \vrule height1.45ex width0.03em depth0em
          \kern-0.30em \rm Q}}
\def\qu{{\mathchoice
    {\sQQ}
    {\sQQ}
  {\kern 0.225em \vrule height1.05ex width0.025em depth0em \kern-0.25em \rm Q}
  {\kern 0.180em \vrule height0.78ex width0.020em depth0em \kern-0.20em \rm Q}
        }}
\def\sCC{{\kern 0.27em \vrule height1.45ex width0.03em depth0em
          \kern-0.30em \rm C}}
\def\complex{{\mathchoice
    {\sCC}
    {\sCC}
  {\kern 0.225em \vrule height1.05ex width0.025em depth0em \kern-0.25em \rm C}
  {\kern 0.180em \vrule height0.78ex width0.020em depth0em \kern-0.20em \rm C}
        }}

\newcommand{\R}{\mathbb{R}}

\newcommand{\ba}{\begin{array}}
\newcommand{\ea}{\end{array}}
\newcommand{\be}{\begin{equation}}
\newcommand{\ee}{\end{equation}}
\newcommand{\bea}{\begin{eqnarray}}
\newcommand{\eea}{\end{eqnarray}}
\newcommand{\beaa}{\begin{eqnarray*}}
\newcommand{\eeaa}{\end{eqnarray*}}

\def\z{\zeta}

\font\tenmath=msbm10 \font\sevenmath=msbm7 \font\fivemath=msbm5
\newfam\mathfam \textfont\mathfam=\tenmath
\scriptfont\mathfam=\sevenmath \scriptscriptfont\mathfam=\fivemath
\def\math{\fam\mathfam}

\def \={{\buildrel {\rm (law)} \over =}}

\def \R{\mathbb{R}}

\def \Q{{\math Q}}

\def\HH{\EuFrak H}
\def\qed{ \hfill \vrule width.25cm height.25cm depth0cm\smallskip}

\newcommand{\cont}[1]{\stackrel{#1}{\frown}}

\newcommand{\basa}{\begin{assumption}}
\newcommand{\easa}{\end{assumption}}

\newcommand{\bas}{\begin{assum}}
\newcommand{\eas}{\end{assum}}

\newcommand{\ignore}[1]{}
\usepackage{tikz}

\begin{document}
\renewcommand{\thefootnote}{\fnsymbol{footnote}}

\renewcommand{\thefootnote}{\fnsymbol{footnote}}
\begin{center}
{\bf SEMICIRCULAR LIMITS ON THE FREE POISSON CHAOS:} \\ 
{{\bf COUNTEREXAMPLES TO A TRANSFER PRINCIPLE}}

\medskip

\normalsize
by Solesne Bourguin\footnote{Universit\'e du Luxembourg. Facult\'e des Sciences, de la Technologie et de la Communication: Unit\'e de Recherche en Math\'ematiques. 6, rue Richard Coudenhove-Kalergi, L-1359 Luxembourg. Email: {\tt solesne.bourguin@gmail.com}} and Giovanni Peccati\footnote{Universit\'e du Luxembourg. Facult\'e des Sciences, de la Technologie et de la Communication: Unit\'e de Recherche en Math\'ematiques. 6, rue Richard Coudenhove-Kalergi, L-1359 Luxembourg. Email: {\tt giovanni.peccati@gmail.com}} \\ {\it  Universit\'e du Luxembourg}\\~\\

\end{center}

{\small \noindent {\bf Abstract}: We establish a class of sufficient conditions, ensuring that a sequence of multiple integrals with respect to a free Poisson measure converges to a semicircular limit. We use this result to construct a set of explicit counterexamples, showing that the \textit{transfer principle} between classical and free Brownian motions (recently proved by Kemp, Nourdin, Peccati and Speicher (2012)) does not extend to the framework of Poisson measures. Our counterexamples implicitly use kernels appearing in the classical theory of random geometric graphs. Several new results of independent interest are obtained as necessary steps in our analysis, in particular: (i) a multiplication formula for free Poisson multiple integrals, (ii) diagram formulae and spectral bounds for these objects, and (iii) a counterexample to the general universality of the Gaussian Wiener chaos in a classical setting.    \\


    
\section{Introduction}

Let $W =\{W_t : t\geq 0\}$ be a real--valued standard Brownian motion defined on the probability space $(\Omega, \mathcal{F}, P)$, and let $S=\{S_t : t\geq 0\}$ be a free Brownian motion defined on the non--commutative probability space $(\mathscr{A}, \varphi)$. Given $q\geq 1$ and a kernel $f\in L^2(\R_+^q, \mu^q)$ (where $\mu$ stands for the Lebesgue measure), we shall denote by $I_q^S(f)$ and $I_q^W(f)$, respectively, the {\it multiple Wigner integral} of $f$ with respect to $S$ (see \cite{bianespeicher}), and the {\it multiple Wiener integral} of $f$ with respect to $W$ (see \cite{np-book}). Objects of this type play crucial roles in non--commutative and classical stochastic analysis; in particular, they constitute the basic building blocks that allow one to define Malliavin operators in both contexts (see again \cite{bianespeicher, np-book}).

\medskip

\noindent The following asymptotic result, that contains both a semicircular limit theorem and a transfer principle, was first proved in \cite{knps}. We denote by $\mathcal{S}(0,1)$ a centered free semicircular random variable with parameter 1, and by $\mathcal{N}(0,1)$ a centered Gaussian random variable with unit variance (all the notions evoked in the introduction will be formally defined in Section \ref{ss:pre} and Section 3).

\begin{theorem}[See \cite{knps}] \label{t:knps} Fix an integer $q\geq 2$, and let $\{f_n : n\geq 1\} \subset L^2(\R_+^q)$ be a sequence of mirror symmetric kernels such that $\| f_n\|_{L^2(\R_+^q)}\to 1$. \begin{itemize}

\item[\rm (A)] As $n\to \infty$, one has that $I^S_q(f_n)$ converges in law to $\mathcal{S}(0,1)$ if and only if $\varphi(I^S_q(f_n)^4) \to \varphi(\mathcal{S}(0,1)^4) = 2$.

\item[\rm (B)]{\rm (Transfer principle)} If each $f_n$ is fully symmetric, then $I^S_q(f_n)$ converges in law to $\mathcal{S}(0,1)$ if and only if $I_q^W(f_n)$ converges in law to $\sqrt{q!}\mathcal{N}(0,1)$.

\end{itemize}

\end{theorem}

\noindent Part (A) of the previous statement is indeed a free counterpart of a \textit{fourth moment central limit theorem} by Nualart and Peccati, originally proved in the framework of multiple Wiener integrals with respect to a Gaussian field (see e.g. \cite{np-jfa, npr-jfa}, as well as \cite[Chapter 5]{np-book}). On the other hand, Part (B) establishes a transfer principle according to which, in the presence of symmetric kernels, convergence to the semicircular distribution on the Wigner chaos is exactly equivalent to a chaotic Gaussian limit on the classical Wiener space. Recall that, in free probability, the semicircular distribution plays the same central role as the Gaussian distribution in the classical context. As shown in \cite{DNN, knps}, such a transfer principle allows one to deduce almost immediately free versions of important central limit theorems from classical stochastic analysis, like for instance non--commutative versions of the famous Breuer--Major Theorem for Gaussian--subordinated sequences (see \cite[Chapter 7]{np-book}). 

\medskip

\noindent Since its appearance, Theorem \ref{t:knps} has triggered a number of generalisations: see e.g. \cite{bps} for multidimensional statements, \cite{np} for convergence results towards the free Poisson distribution, \cite{nyet} for alternate proofs, \cite{aria} for fourth moment theorems in the presence of freely infinitely divisible distributions,  \cite{DN} for limit theorems towards the {\it tetilla law}, \cite{DNN} for extensions to the $q$--Brownian setting, and \cite{noupoly} for a characterisation of convergence in law within the second Wigner and Wiener chaos. 

\medskip

\noindent In the present paper, we are interested in the following natural question: {\it do Part {\rm (A)} and {\rm (B)} of Theorem \ref{t:knps} continue to hold whenever one replaces $S$ with a centered free Poisson measure on the real line (see e.g. {\rm \cite{ans2000, ans2002, ans1}}), and $W$ with a centered classical Poisson process?} Our answer will be that, while Part (A) of Theorem \ref{t:knps} can be suitably extended to the free Poisson framework, the transfer principle of Part (B) fails for every order $q$ of integration. 

\medskip

\noindent We denote by $\{\hat{N}_t : t\geq 0\}$ and $\{\hat{\eta}_t : t\geq 0\}$, respectively, a free centered Poisson process, and a classical centered Poisson process. Given $q\geq 1$ and a kernel $f\in L^2(\R_+^q, \mu^q)$, we write $I_q^{\hat{N}}(f)$ and $I_q^{\hat{\eta}}(f)$, respectively, to indicate the free multiple integral of $f$ with respect to $\hat{N}$ (see \cite{ans0, ans1}), and the multiple Wiener integral of $f$ with respect to $\hat{\eta}$ (see \cite{PecTaq}). Formally, our aim in this paper is to prove a general version of the following statement (the generality consists in the fact that we shall consider Poisson measures on $\R^d$):

\begin{theorem} Let $q\geq 2$ be an arbitrary integer. 
\label{t:bp}
\begin{itemize}

\item[\rm (A)] Let $\{f_n : n\geq 1\} \subset L^2(\R_+^q)$ be a { tamed} sequence of mirror symmetric kernels such that $\| f_n\|_{L^2(\R_+^q)}\to 1$,
as $n\to \infty$. Then, $I^{\hat{N}}_q(f_n)$ converges in law to $\mathcal{S}(0,1)$ if and only if $\varphi(I^{\hat{N}}_q(f_n)^4) \to 2$.

\item[\rm (B)]{\rm (Counterexamples to a transfer principle)} There exists a tamed sequence of fully symmetric kernels $\{g_n : n\geq 1\}\subset L^2(\R_+^q)$ such that, as $n\to\infty$, $I^{\hat{N}}_q(g_n)$ converges in distribution to $\mathcal{S}(0,1)$, while $I^{\hat{\eta}}_q(g_n)$ converges in distribution to a centered Poisson random variable.

\end{itemize}

\end{theorem}

\noindent The notion of a {\it tamed sequence} of kernels is introduced in Definition \ref{d:tamed} below: this additional assumption is needed in order to deal with the complicated combinatorial structures arising from the computation of cumulants. As a by-product of our analysis and of the findings of \cite{LacPec}, we shall also prove that the following partial counterpart to Theorem \ref{t:bp}--(B) holds for every order $q$ of integration: assume that the fully symmetric tamed sequence $\{f_n \colon n \geq 1\}$ is composed of kernels with a constant sign; if $I^{\hat{\eta}}_q(f_n)$ converges in law to a normal random variable, one has necessarily that $I^{\hat{N}}_q(f_n)$ converges in law to a semicircular limit.

\medskip

\noindent Part (A) of Theorem \ref{t:bp} is a genuine generalization of the first part of Theorem \ref{t:knps} in a free Poisson framework. It also provides a free counterpart to the central limit theorems for classical multiple Poisson integrals proved in \cite{LacPec, PSTU}: in particular, it is a remarkable fact that -- unlike in the classical framework -- we are able to obtain a genuine \textit{fourth moment theorem} without assuming that the kernels $f_n$ have a constant sign.

\medskip

\noindent In order to show this result, we will overcome a number of technical difficulties: indeed, albeit the foundations of free stochastic integration with respect to a free Poisson measure have been firmly established in the influential papers \cite{ans2000, ans2002, ans0, ans1}, none of these references provides an explicit treatment of product and diagram formulae for general kernels, which are the basic ingredients of our techniques. We devote the whole Section \ref{s:alg} to filling this gap. In particular:
\begin{itemize}

\item[--] In Section \ref{ss:prod} and Section \ref{ss:ext}, we establish general multiplication formulae for free multiple Poisson integrals, based on the use of \textit{star contractions} operators. Such formulae are in some sense analogous to the ones  available in the classical case: nonetheless, some remarkable differences will appear, basically motivated by the rigid combinatorial structure of the lattice of non--crossing partitions. Our proofs are based on the Fock space representation of free random measures.

\item[--] In Section \ref{ss:dia} and Section \ref{ss:ext} we will prove new combinatorial diagram formulae for free multiple Poisson integrals. Such formulae are particularly useful for deducing explicit spectral bounds. Our proofs are based on computations of joint free cumulants (see \cite{NicSpe} for a detailed exposition of this notion), and can be seen as a substantial generalization of the combinatorial cumulant and moment formulae for the centered free Poisson distribution provided in \cite{np}.
\end{itemize}

\noindent The tamed sequence appearing in Part (B) of Theorem \ref{t:bp} is explicitly constructed in Section \ref{s:cont}: we will see that the kernels $g_n$ are intimately related to objects appearing in the classical theory of random graphs.

\medskip

\noindent The rest of the paper is organized as follows: Section \ref{ss:pre} is devoted to preliminary notions of non--commutative stochastic analysis. Section 3 studies the algebra of free multiple integrals. Section 4 contains the statement and proof of our main limit theorems, whereas Section 5 provides an explicit construction of the examples and counterexamples mentioned above. Section 6 contains the (quite technical) proofs of product and moment/cumulant formulae.

\section{Elements of free probability}\label{ss:pre}
In this section, we recall some basic notions of free probability. The reader is referred to the two fundamental references \cite{NicSpe, vdn} for any unexplained definition or result.
\\~\\
{\bf Non--commutative probability spaces.} A {\it non--commutative probability space} is a pair $(\mathscr{A}, \varphi)$, where $\mathscr{A}$ is a unital $\ast$--algebra, and $\varphi : \mathscr{A}\to \mathbb{C}$ is a unitary linear operator, which is moreover positive (meaning that $\varphi(YY^*) \geq 0$ for every $Y \in \mathscr{A}$). We shall refer to the self--adjoint elements of $\mathscr{A}$ as {\it random variables}.  In what follows, we will sometimes evoke the more specialized notions of a $C^*$-- and $W^*$--{\it probability space}. A $C^*$--probability space is a non--commutative probability space such that $\mathscr{A}$ is a $C^*$--algebra. A $W^*$--probability space is a non--commutative probability space such that $\mathscr{A}$ is a von--Neumann algebra,
and $\varphi : \mathscr{A}\to \mathbb{C}$ is a linear operator satisfying the following assumptions: $\varphi$ is weakly continuous, unital, faithful (meaning that $\varphi(YY^\ast) = 0 \Rightarrow Y = 0$, for every $Y\in \mathscr{A}$), positive and tracial (meaning that $\varphi(XY) = \varphi(YX)$ for all $X,Y \in \mathscr{A}$). Whenever the trace property is satisfied, one customarily says that $\varphi$ is a {\it tracial state} on $\mathscr{A}$. 
\\~\\
\noindent{\bf Free distributions and convergence.}
Let $(\mathscr{A},\varphi)$ be a non--commutative probability space and let $X\in \mathscr{A}$. The $k^{\mbox{\tiny{th}}}$ {\it moment} of $X$ is given by the quantity $\varphi(X^{k})$, $k=0,1,...$. Now assume that $X$ is a self--adjoint bounded element of $\mathscr{A}$ (in other words, $X$ is a bounded free random variable), and write $\rho(X)\in [0, \infty)$ to indicate the {\it spectral radius} of $X$. The {\it law} (or {\it spectral measure}) of $X$ is defined as the unique Borel probability measure $\mu_{X}$ on the real line such that: (i) $\mu_X$ has support in $[-\rho(X), \rho(X)]$, and (ii)
$\int_{\mathbb{R}}P(t)\ d\mu_{X}(t) = \varphi(P(X))$
for every polynomial $P \in \mathbb{R}\left[ X\right]$. The existence and uniqueness of $\mu_X$ in such a general framework are proved e.g. in \cite[Theorem 2.5.8]{tao} (see also \cite[Proposition 3.13]{NicSpe}). We shall sometimes use the notation $X \sim \mu_X$ to indicate that $\mu_X$ is the law of $X$. Note that, since $\mu_X$ has compact support, the measure $\mu_X$ is completely determined by the sequence $\left\lbrace \varphi(X^k) \colon k\geq 1\right\rbrace $. 
\\~\\
\noindent{\bf Convergence.} Let $\left\lbrace X_{n} \colon n \geq 1\right\rbrace $ be a sequence of non--commutative random variables, each possibly belonging to a different non--commutative probability space $(\mathscr{A}_n, \varphi_n)$. The sequence $\left\lbrace X_{n} \colon n \geq 1\right\rbrace $ is said to {\it converge in law} to a limiting non--commutative random variable $X_{\infty}$ (defined on $(\mathscr{A}_\infty, \varphi_\infty)$) if $\lim_{n \to+\infty}\varphi_n(P(X_{n})) = \varphi_\infty(P(X_{\infty}))$ for every polynomial $P\in \R[X]$. If $X_\infty, X_n$, $n\geq 1$, are bounded (and therefore the spectral measures $\mu_{X_n}, \mu_{X_\infty}$ are well--defined) this last relation is equivalent to saying that $\int_\R P(t)\, \mu_{X_n}(dt) \to \int_\R P(t)\, \mu_{X_\infty}(dt)$; an application of the method of moments (see e.g. \cite[p. 412]{Bill}) yields immediately that, in this case, one has also that $\mu_{X_n}$ weakly converges to $\mu_{X_\infty}$, that is: $\mu_{X_n}(f) \to \mu_{X_\infty}(f)$, for every $f: \R\to \R$ bounded and continuous. We will sometimes indicate that $X_n$ converges in law to $X_\infty$ by writing $X_{n} \longrightarrow X_{\infty}$ (whenever $X_\infty, X_n$, are bounded, we shall equivalently write $\mu_{X_n} \longrightarrow \mu_{X_\infty}$). Given a non--commutative probability space $(\mathscr{A}, \varphi)$, we write $L^2(\varphi) :=L^2(\mathscr{A},\varphi)$ to indicate the $L^2$ space obtained as the completion of  $\mathscr{A}$ with respect to the norm $\| a\|_2 = \varphi(a^*a )^{1/2}$. The following elementary result is used in several occasions throughout the paper. 
\begin{lemma}\label{l:lot} Let $\{ F,F_n : n\geq 1\}$ be a collection of bounded random variables on a tracial non--commutative probability space $(\mathscr{A}, \varphi)$, and assume that $F_n \to F$ in $L^2(\varphi)$ and $R:= \sup_n\rho(F_n) <\infty$. Then, one has that {\rm (i)} $F^m_n \to F^m$ in $L^2(\varphi)$, and {\rm (ii)} $\varphi(F_n^m) \to \varphi(F^m)$, for every integer $m\geq 1$.
\end{lemma}
\noindent{\bf Proof.} Since $| \varphi(F_n^m)- \varphi(F^m)| \leq \|F^m_n-F^m\|_2$ (by Cauchy--Schwarz), it suffices to prove point (i). By assumption, the conclusion is true for $m=1$. We shall apply a recursive argument, and show that, if relation (i) holds for $m$, then it also holds for $m+1$. To prove this, we use the following identity, satisfied by any pair of non--commutative variables $x,y$,
\begin{eqnarray*}
 x^{m+1} - y^{m+1}  &=& \frac14 \big\{ (x^{m}+y^{m})(x-y) + (x^{m}-y^{m})(x+y) \\
 &&\quad\quad\quad + (x-y)(x^{m}+y^{m})+(x+y)(x^{m}-y^{m})\big\}.
\end{eqnarray*}
We apply the last relation to $x= F_n$  and $y=F$ to rewrite one the quantity $(F^m_n-F^m)$ appearing in the expression $\|F_n^{m+1} - F^{m+1}\|_2^2 = \varphi((F^m_n-F^m)(F^m_n-F^m))$. Using the tracial property of $\varphi$ and exploiting several times the Cauchy--Schwarz inequality together with the fact that $R<\infty$, one sees eventually that
\[
\|F_n^{m+1} - F^{m+1}\|_2^2 = O(1) \big\{ \|F_n - F\|_2+\|F_n^{m} - F^{m}\|_2\big\},
\]
where $O(1)$ stands for a positive bounded numerical sequence (possibly depending on $m$). This concludes the proof.
\qed
\\~\\
{\bf Free independence.}
Another central concept in non--commutative probability is that of \textit{free independence}. It is the counterpart of independence in a classical probability setting. Consider a free tracial probability space $(\mathscr{A},\varphi)$, and let $\mathscr{A}_{1},\ldots,\mathscr{A}_{p}$ be unital subalgebras of $\mathscr{A}$. The family of subalgebras $\left\lbrace \mathscr{A}_{1},\ldots,\mathscr{A}_{p} \right\rbrace $ is called {\it free} or {\it freely independent} if $\varphi(X_{1}X_{2}\cdots X_{m}) = 0$, whenever the following conditions are met:
the $X_{1},\ldots,X_{m}$ are such that, for each $1\leq j\leq m$, $X_{j} \in \mathscr{A}_{i(j)}$, where $i(1) \neq i(2)$, $i(2) \neq i(3)$, $\ldots$, $i(m-1) \neq i(m)$, and  $\varphi(X_{j}) = 0$ for all $j$. Random variables are termed free or freely independent if the unital algebras they generate are free. The concept of free independence deviates considerably from the notion of classical independence: for example, if $X$ and $Y$ are free, we have $\varphi(XYXY) = \varphi(Y)^{2}\varphi(X^{2}) + \varphi(X)^{2}\varphi(Y^{2}) - \varphi(X)^{2}\varphi(Y)^{2}$, which is in contrast to the relation $\E[XYXY]=\E[X^2]\E[Y^2]$, holding for independent (and square--integrable) random variables defined on a classical probability space $(\Omega, \mathscr{F}, \mathbb{P})$.
\\~\\
{\bf Non--crossing partitions.}  Given an integer $m \geq 1$, we write $[m] = \{1,\ldots,m\}$. A {\it partition} of $[m]$ is a collection of non--empty and disjoint subsets of $[m]$, called {\it blocks}, such that their union is equal to $[m]$. The cardinality of a block is called {\it size}. A block is said to be a {\it singleton} if it has size one. We adopt the convention of ordering the blocks of a given partition $\pi = \{B_1,\ldots,B_r\}$ by their least element, that is: $\min\, B_i < \min\, B_j$ if and only if $i< j$. A partition $\pi$ of $[n]$ is said to be {\it non--crossing} if one cannot find integers $p_1,q_1,p_2,q_2$ such that: (a) $1\leq p_1 < q_1 < p_2 < q_2\leq m$, (b) $p_1,p_2$ are in the same block of $\pi$, (c)  $q_1,q_2$ are in the same block of $\pi$, and (d) the $p_i$'s are not in the same block of $\pi$ as the $q_i$'s. The collection of the non--crossing partitions of $[n]$ is denoted by $NC(n)$, $n\geq 1$. Fig. 1 provides a standard graphical representation of the two partitions of $[10]$ given by $\pi_1 = \{\{1,9\}, \{2,10\}, \{3,4\}, \{5,8\}, \{6,7\}\}$ (crossing) and $\pi_2 = \{\{1,10\}, \{2,9\}, \{3,4\}, \{5,8\}, \{6,7\}\}$ (non--crossing). In particular, the crossing nature of $\pi_1$ is determined by the two blocks $\{1,10\}$ and $\{2,9\}$ (colored in light blue). 
\begin{figure}[H]
\begin{center}
\label{f:b}
\scalebox{1.5} 
{
\begin{pspicture}(0,-0.43661112)(8.904111,0.4425)
\definecolor{color2411}{rgb}{0.6,0.8,1.0}
\definecolor{color1842}{rgb}{0.00392156862745098,0.00392156862745098,0.00392156862745098}
\psdots[dotsize=0.09](2.445,-0.3775)
\psdots[dotsize=0.09,linecolor=color1842](2.445,-0.3775)
\psdots[dotsize=0.09,linecolor=color1842](3.245,-0.3775)
\psdots[dotsize=0.09,linecolor=color1842](1.245,-0.3775)
\psline[linewidth=0.04cm,linecolor=color2411](0.045,0.4225)(3.245,0.4225)
\psline[linewidth=0.04cm,linecolor=color2411](3.245,0.4225)(3.245,-0.3775)
\psdots[dotsize=0.09](0.045,-0.3775)
\psdots[dotsize=0.09](0.445,-0.3775)
\psdots[dotsize=0.09](0.845,-0.3775)
\psdots[dotsize=0.09](1.645,-0.3775)
\psdots[dotsize=0.09](2.045,-0.3775)
\psdots[dotsize=0.09](2.845,-0.3775)
\psdots[dotsize=0.09](3.645,-0.3775)
\psline[linewidth=0.04cm,linecolor=color2411](0.045,0.4225)(0.045,-0.3775)
\psline[linewidth=0.04cm](1.645,0.0225)(2.845,0.0225)
\psline[linewidth=0.04cm](2.845,0.0225)(2.845,-0.3775)
\psline[linewidth=0.04cm](1.645,0.0225)(1.645,-0.3775)
\psline[linewidth=0.04cm](0.845,-0.1775)(0.845,-0.3775)
\psline[linewidth=0.04cm](0.845,-0.1775)(1.245,-0.1775)
\psline[linewidth=0.04cm](1.245,-0.1775)(1.245,-0.3775)
\psline[linewidth=0.04cm](2.045,-0.3775)(2.045,-0.1775)
\psline[linewidth=0.04cm](2.045,-0.1775)(2.445,-0.1775)
\psline[linewidth=0.04cm](2.445,-0.1775)(2.445,-0.3775)
\psline[linewidth=0.04cm,linecolor=color2411](0.445,0.2225)(3.645,0.2225)
\psline[linewidth=0.04cm,linecolor=color2411](0.445,0.2225)(0.445,-0.3775)
\psline[linewidth=0.04cm,linecolor=color2411](3.645,0.2225)(3.645,-0.3775)
\psdots[dotsize=0.09](0.045,-0.3775)
\psdots[dotsize=0.09](0.445,-0.3775)
\psdots[dotsize=0.09](3.245,-0.3775)
\psdots[dotsize=0.09](3.645,-0.3775)
\psdots[dotsize=0.09](0.845,-0.3775)
\psdots[dotsize=0.09](1.245,-0.3775)
\psdots[dotsize=0.09](1.645,-0.3775)
\psdots[dotsize=0.09](2.045,-0.3775)
\psdots[dotsize=0.09](2.445,-0.3775)
\psdots[dotsize=0.09](2.845,-0.3775)
\psdots[dotsize=0.09](7.645,-0.3775)
\psdots[dotsize=0.09,linecolor=color1842](7.645,-0.3775)
\psdots[dotsize=0.09,linecolor=color1842](8.445,-0.3775)
\psdots[dotsize=0.09,linecolor=color1842](6.445,-0.3775)
\psline[linewidth=0.04cm,linecolor=color2411](5.245,0.4225)(8.845,0.4225)
\psline[linewidth=0.04cm,linecolor=color2411](8.845,0.4225)(8.845,-0.3775)
\psdots[dotsize=0.09](5.245,-0.3775)
\psdots[dotsize=0.09](5.645,-0.3775)
\psdots[dotsize=0.09](6.045,-0.3775)
\psdots[dotsize=0.09](6.845,-0.3775)
\psdots[dotsize=0.09](7.245,-0.3775)
\psdots[dotsize=0.09](8.045,-0.3775)
\psdots[dotsize=0.09](8.845,-0.3775)
\psline[linewidth=0.04cm,linecolor=color2411](5.245,0.4225)(5.245,-0.3775)
\psline[linewidth=0.04cm](6.845,0.0225)(8.045,0.0225)
\psline[linewidth=0.04cm](8.045,0.0225)(8.045,-0.3775)
\psline[linewidth=0.04cm](6.845,0.0225)(6.845,-0.3775)
\psline[linewidth=0.04cm](6.045,-0.1775)(6.045,-0.3775)
\psline[linewidth=0.04cm](6.045,-0.1775)(6.445,-0.1775)
\psline[linewidth=0.04cm](6.445,-0.1775)(6.445,-0.3775)
\psline[linewidth=0.04cm](7.245,-0.3775)(7.245,-0.1775)
\psline[linewidth=0.04cm](7.245,-0.1775)(7.645,-0.1775)
\psline[linewidth=0.04cm](7.645,-0.1775)(7.645,-0.3775)
\psline[linewidth=0.04cm,linecolor=color2411](5.645,0.2225)(8.445,0.2225)
\psline[linewidth=0.04cm,linecolor=color2411](5.645,0.2225)(5.645,-0.3775)
\psline[linewidth=0.04cm,linecolor=color2411](8.445,0.2225)(8.445,-0.3775)
\psdots[dotsize=0.09](5.245,-0.3775)
\psdots[dotsize=0.09](5.645,-0.3775)
\psdots[dotsize=0.09](8.445,-0.3775)
\psdots[dotsize=0.09](8.845,-0.3775)
\psdots[dotsize=0.09](6.045,-0.3775)
\psdots[dotsize=0.09](6.445,-0.3775)
\psdots[dotsize=0.09](6.845,-0.3775)
\psdots[dotsize=0.09](7.245,-0.3775)
\psdots[dotsize=0.09](7.645,-0.3775)
\psdots[dotsize=0.09](8.045,-0.3775)
\end{pspicture} 
}
\caption{\sl A crossing (left) and a non--crossing (right) partition of $[10]$.}
\end{center}
\end{figure}
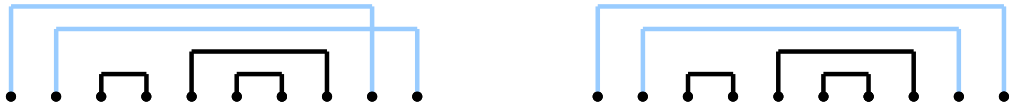

\noindent{\bf Lattice structure of $NC(n)$.} It is a well--known fact (see e.g. \cite[p. 144]{NicSpe}) that the reversed refinement order (written $\preceq$) induces a lattice structure on $NC(n)$: we shall denote by $\vee$ and $\wedge$, respectively, the associated {\it join} and {\it meet} operations, whereas $\hat 0 = \{\{1\},\ldots,\{n\}\}$ and $\hat 1=\{[n]\}$ are the corresponding minimal and maximal partitions of the lattice. Note that the meet operation $\wedge$ in $NC(n)$ coincides with the meet operation in the larger lattice $\mathcal{P}(n)$ of all partitions of $[n]$, but the same property does not hold for the join operation (see \cite[Remark 9.19]{NicSpe}). However, if $\pi^*$ is a block partition such as the ones defined in Definition \ref{respectfulpartitions} below and $\sigma\in NC(n)$, then the join $\pi^*\vee \sigma$ is the same as the one obtained by regarding $\pi^*, \sigma$ as elements of $\mathcal{P}(n)$ (see \cite[Exercise 9.43]{NicSpe}). We record here a useful characterization of non--crossing partitions (see \cite[Remark 9.2]{NicSpe}). 
\begin{lemma}\label{l:1}
A partition $\pi$ of $\left[n \right]$ is non--crossing if and only if at least one block $B \in \pi$ is an interval and $\pi\setminus B$ is non--crossing, that is: $B \in \pi$ has the form $B = \left\lbrace k, k+1,\ldots,k+p \right\rbrace $ for some $1 \leq k \leq n$ and $p \geq 0$ (with $k+p \leq n$) and one has the relation $$\pi \setminus B \in {NC}\left( \left\lbrace 1,\ldots,k-1,k+p+1,\ldots,n\right\rbrace \right) \cong {NC}\left( n -(p+1) \right).$$
\end{lemma}
The last relation in the statement simply means that the partition of $$\gamma := \left\lbrace 1,\ldots,k-1,k+p+1,\ldots,n\right\rbrace$$ obtained by deleting $B$ from $\pi$ is non--crossing, and that the class of non--crossing partitions of the set $\gamma$ constitutes a lattice isomorphic to ${NC}\left( n -(p+1) \right)$.
\\~\\
{\bf Catalan numbers.}  For every $n\geq 1$, the quantity $C_n = | NC(n) |$, where $ |A |$ indicates the cardinality of a given set $A$, is called the $n^{\mbox{\tiny{th}}}$ {\it Catalan number}. One sets by convention $C_0 = 1$. Also, recall the explicit expression $C_n = \frac{1}{n+1}\binom{2n}{n}$, from which one infers that
\begin{equation}\label{e:stir}
C_n \sim \pi^{-1/2} n^{-3/2}4^n, \quad n\to \infty,
\end{equation}
where we have made a standard use of Stirling's formula.
\\~\\
{\bf Free cumulants.} Given a random variable $X$ in a non--commutative probability space $(\mathscr{A}, \varphi)$, we denote by $\{\kappa_m(X) : m\geq 1\}$ the sequence of the {\it free cumulants}  of $X$. We recall 
(see \cite[p. 175]{NicSpe}) that the free cumulants of $X$ are defined by the following relation: for every $m\geq 1$,
\begin{equation}\label{e:momcum}
\varphi(X^m) = \sum_{\pi = \{B_1,\ldots,B_r\}\in NC(m)} \,\,\prod_{i=1}^r\kappa_{|B_i|}(X),
\end{equation}
where $|B_i|$ indicates the size of the block $B_i$ of the non--crossing partition $\pi$. It is clear from (\ref{e:momcum}) that the sequence  $\{\kappa_m(X) : m\geq 1\}$ completely determines the moments of $X$ (and viceversa).
\\~\\
{\bf Semicircular distribution.}
The centered {\it semicircular distribution} of parameter $t>0$, denoted by $\mathcal{S}(0,t)(dx)$, is the probability distribution given by $$\mathcal{S}(0,t)(dx) = (2\pi t)^{-1}\sqrt{4t-x^2}\,dx, \quad |x|< 2\sqrt{t}.$$ We recall the classical relation: $$\int_{-2\sqrt{t}}^{2\sqrt{t}} x^{2m} \mathcal{S}(0,t)(dx) = C_m t^m,$$ where $C_m$ is the $m^{\mbox{\tiny{th}}}$ Catalan number (so that e.g. the second moment of $\mathcal{S}(0,t)$ is $t$). Since the odd moments of $\mathcal{S}(0,t)$ are all zero, one deduces from the previous relation and (\ref{e:momcum}) (e.g. by recursion) that the free cumulants of a random variable $s$ with law $\mathcal{S}(0,t)$ are all zero, except for $\kappa_2(s) =\varphi(s^2)= t$. We recall (see e.g. \cite{ans1}) that the class of orthogonal polynomials associated with $\mathcal{S}(0,t)$ is given by the (generalized) {\it Tchebycheff polynomials} $\{H_{0,m} (x,t) : m=0,1,\ldots\}$, defined by the recurrence relations 
\begin{numcases}
~H_{0,0}(x,t) = 1, \nonumber \\
H_{0,1}(x,t) = x, \nonumber \\
xH_{0,m}(x,t) = H_{0,m+1}(x,t)  + tH_{0,m-1}(x,t). \nonumber
\end{numcases}
\\~\\
{\bf Free Poisson distribution.} The {\it free Poisson distribution} with rate $\lambda >0 $, denoted by $P(\lambda)$, is the probability distribution defined 
as follows: (a) if $\lambda\in (0,1]$, then $P(\lambda) = (1-\lambda)\delta_0 + \lambda\widetilde{\nu}$, and (b) if $\lambda > 1$, then 
$P(\lambda) = \widetilde{\nu}$, where $\delta_0$ stands for the Dirac mass at $0$. Here, $\widetilde{\nu}(dx) = 
(2\pi x)^{-1}\sqrt{4\lambda-(x-1-\lambda)^2}dx, \, x\in \big((1-\sqrt{\lambda})^2,(1+\sqrt{\lambda})^2\big)$. If $X_\lambda$ has the $P(\lambda)$ distribution, then \cite[Proposition 12.11]{NicSpe} implies that 
\begin{equation}\label{e:cumPOISS}
\kappa_m(X_\lambda) = \lambda, \quad m\geq 1.
\end{equation}
We will denote by $\hat{P}(\lambda)$ the {\it centered Poisson distribution} of parameter $\lambda$, that is: $\hat{P}(\lambda)$ is the law of the free random variable $Z_\lambda := X_\lambda - \lambda 1$, where $1$ is the unit of $\mathscr{A}$. Plainly, $\kappa_1(Z_\lambda) = \varphi(Z_\lambda) = 0$, and $\kappa_m(Z_\lambda) = \lambda = \kappa_m(X_\lambda)$ for every $m\geq 2$. See \cite[Proposition 2.4]{np} for a characterization of the moments of $\hat{P}(\lambda)$ in terms of {\it Riordan numbers}. We also recall (see again \cite{ans1}) that the class of orthogonal polynomials associated with $\hat{P}(\lambda)$ is given by the (generalized) {\it centered Charlier polynomials} $\{C_{0,m} (x,\lambda) : m=0,1,\ldots\}$, defined by the recurrence relations 
\begin{numcases}
~C_{0,0}(x,\lambda) = 1, \nonumber \\
C_{0,1}(x,\lambda) = x, \nonumber \\
xC_{0,m}(x,\lambda) = C_{0,m+1}(x,\lambda)  +C_{0,m}(x,\lambda)+  \lambda \,C_{0,m-1}(x,\lambda). \nonumber
\end{numcases}
\\~\\
{\bf Free Brownian motion.}
A free Brownian motion $\left\lbrace S(t) \colon t\geq 0\right\rbrace $ is a non--commutative stochastic process. 
It is a family of self--adjoint operators defined on a non--commutative probability space, having the following characteristic properties:
\begin{enumerate}
\item $S(0) = 0$.
\item For $0\leq t_1\leq t_2$, the law of $S(t_{2}) - S(t_{1})$ is the semicircular distribution ${\cal S}(0, t_2-t_1)$.
\item For all $n$ and $0<t_{1}< t_{2}<\cdots <t_{n}$, 
the increments $S(t_{1})$, $S(t_{2}) - S(t_{1})$, $\ldots$, $S(t_{n}) - S(t_{n-1})$ are freely independent.
\end{enumerate}
~\\
{\bf Free random measures.} Let $(Z,\mathcal{Z})$ be a Polish space, with $\mathcal{Z}$ the associated Borel $\sigma$--field, and let $\mu$ be a positive $\sigma$--finite measure over $(Z,\mathcal{Z})$ with no atoms. We denote by $\mathcal{Z}_{\mu}$ the class of those $A \in \mathcal{Z}$ such that $\mu(A) < \infty$. Let $(\mathscr{A},\varphi)$ be a free tracial probability space and let $\mathscr{A}_{+}$ denote the cone of positive operators in $\mathscr{A}$. Then, a \textit{free semicircular random measure} (resp. a \textit{free Poisson random measure}) with control $\mu$ on $(Z,\mathcal{Z})$ with values in $(\mathscr{A},\varphi)$ is a mapping $S : \mathcal{Z}_{\mu} \rightarrow \mathscr{A}$ (resp. $N : \mathcal{Z}_{\mu} \rightarrow \mathscr{A}_{+}$), with the following properties:
\begin{enumerate}
\item For any set $A$ in $\mathcal{Z}_{\mu}$, $S(B)$ (resp. $N(B)$) is a centered semicircular (resp. a free Poisson) random variable with variance (resp. parameter) $\mu(A)$;
\item If $r \in \mathbb{N}$ and $A_1,\ldots,A_r \in \mathcal{Z}_{\mu}$ are disjoint, then $S(A_1),\ldots,S(A_r)$ (resp. $N(A_1),\ldots,N(A_r)$) are free;
\item If $r \in \mathbb{N}$ and $A_1,\ldots,A_r \in \mathcal{Z}_{\mu}$ are disjoint, then $S\left( \bigcup_{j=1}^{r}A_j\right) = \sum_{i=1}^{r}S\left( A_j\right)$ (resp. $N\left( \bigcup_{j=1}^{r}A_j\right) = \sum_{i=1}^{r}N\left( A_j\right)$).
\end{enumerate}
The existence of free Poisson measures on an appropriate $W^*$--probability space is guaranteed by \cite[Theorem 3.3]{BT} and \cite[Theorem 5.1]{BAV} (see also \cite[Chapter 4]{kar}). The existence of a free semicircular measure on a tracial faithful $C^*$--probability space follows e.g. from the discussion contained in \cite[pp. 102--108]{NicSpe}. Note that a free Brownian motion can be identified with a semicircular measure on $(\R_+, \mathscr{B}(\R_+))$ and control given by the Lebesgue measure. If $N$ is a free Poisson measure with control $\mu$, we will denote by $\hat N$ the collection of random variables
$$
\hat{N}(A) = N(A) - \mu(A)1, \quad A\in \mathcal{Z}_\mu,
$$
where $1$ is the unit of $\mathscr{A}$, in such a way that $\hat{N}(A)\sim \hat{P}(\mu(A))$. We call $\hat N$ a {\it free centered Poisson measure}.

\section{The algebra of free Poisson multiple integrals}\label{s:alg}

We will now provide a detailed description of the algebraic structure of the class of multiple integrals (of bounded kernels with bounded support) with respect to a free Poisson measure. Since the present paper is meant as a partial counterpart to the findings of \cite{DN, DNN, knps, np}, we decide to attach to our discussion a parallel (succinct) analysis of multiple integration with respect to a free Brownian motion: in particular, this choice allows one to better appreciate some important combinatorial differences between the non--commutative Poisson and Brownian settings. The proofs of the main results of this section are detailed in Section \ref{s:proofs}.

\subsection{Further notation}

For definiteness, for the rest of the paper we shall focus on the case where $$(Z,\mathcal{Z}) = \left(\R^d, \mathscr{B}(\R^d)\right),\quad \mbox{for some } d\geq1,$$ and $\mu$ equals the $d$--dimensional Lebesgue measure; we will write $L^2(Z^1) = L^2(Z) := L^2(Z,\mathcal{Z}, \mu)$ for the space of complex--valued functions that are square--integrable with respect to $\mu$.  The following additional standard notational conventions are in order:

\begin{itemize}

\item[\rm (1)] For every integer $q\geq 2$, the space $L^2(Z^q)$ is the collection of all complex--valued functions on $Z^q$ that are square--integrable with respect to $\mu^q$. Given $f\in L^2(Z^q)$, we write 
\[
f^*(t_1,t_2,\ldots,t_q) = \overline{f(t_q,\ldots,t_2,t_1)},
\]  
and we call $f^*$ the {\it adjoint} of $f$. We say that an element of $L^2(Z^q)$ is {\it mirror symmetric} if
\[
f(t_1,\ldots,t_q) = f^*(t_1,\ldots,t_q),
\]
for almost every vector $(t_1,\ldots,t_q)\in Z^q$. Observe that mirror symmetric functions constitute a Hilbert subspace of $L^2(Z^q)$.
\item[\rm(2)] For every $q\geq 1$ we denote by $\mathscr{E}_q$ the class of {\it elementary elements} of $L^2(Z^q)$, defined as the complex vector space generated by functions of the type 
\begin{equation}\label{e:passa}
f = f_1\otimes \cdots \otimes f_q= \mathds{1}_{A_1}^{\otimes k_1} \otimes \mathds{1}_{A_2}^{\otimes k_2}\otimes \cdots \otimes \mathds{1}_{A_l}^{\otimes k_l},
\end{equation}
where $k_1+\cdots +k_l = q$, each set $A_j$ is bounded (i.e., it is contained in some open ball centered at the origin), and $A_j\cap A_{j+1} = \emptyset$ for $j=1,\ldots,l-1$. We shall also denote by $\mathscr{E}_q^0$ the subset of $\mathscr{E}_q$ composed of the linear combinations of those kernels of the type \eqref{e:passa} such that $l=q$ (and therefore $k_1=\cdots =k_q = 1$). Exploiting the fact that $\mu$ has no atoms, it is clear that $\mathscr{E}_q^0$ (and consequently $\mathscr{E}_q$) is a dense subspace of $L^2(Z^q)$. In several occasions, we will also implicitly use the following fact, which follows from standard properties of the Lebesgue measure.
\begin{lemma}\label{l:taborn} For every bounded measurable $f : Z^q \to \mathbb{C}$ with bounded support, there exists a sequence $\{f_n\} \subset \mathscr{E}_q^0$ such that: {\rm (a)} the support of each $f_n$ is contained in the support of $f$, {\rm (b)} $|f_n|\leq |f|$, and {\rm (c)} $f_n(x)$ converges to $f(x)$ for $\mu$--almost every $x$, and {\rm (d)} $f_n\to f$ $\!\!$ in $L^2(Z^q)$ (by dominated convergence). Moreover, since $\mu$ is non--atomic (and therefore $\mu^q$ does not charge diagonals), one can take each $f_n$ to be an element of $\mathscr{E}_q^{00}$, where $\mathscr{E}_q^{00}$ is defined as the span of those elements of $\mathscr{E}_q^0$ having the `purely non--diagonal' form
\[
\mathds{1}_{A_1}\otimes \cdots \otimes \mathds{1}_{A_q},
\]
where $A_i\cap A_k = \emptyset$ for $i\neq k$. 
\end{lemma}
To keep our notation as compact as possible, we also use the convention $\mathscr{E}_0 =\mathscr{E}^0_0=\mathscr{E}^{00}_0 =\mathbb{C}$.
\item[\rm(4)] Let $f \in L^2\left(Z^m\right)$ and $g \in L^2\left(Z^n\right)$. The general definition of the \textit{star contraction} between $f$ and $g$ goes as follows: for $1 \leq k \leq m \wedge n$ and $j\in \{0,\ldots,k\}$, we set 
\begin{eqnarray*}
&& f  \star_{k}^{j} g (t_1,\ldots, t_{m+n - 2k+j}) = \\
&& \quad \int_{Z^{k-j}}f(t_1,\ldots ,t_{m-k+1},\ldots, t_{m-k+j}, s_{k-j}, \ldots , s_1) \\
&& \qquad\quad\quad \times g(s_1, \ldots , s_{k-j} , t_{m-k+1},\ldots,t_{m-k+j},\ldots,  t_{m+n - 2k+j})\mu(ds_1)\cdots \mu (ds_{k-j}).
\end{eqnarray*}
\noindent This general notation (that will be used e.g. in Section 5) corresponds to the one adopted in classical stochastic analysis (see \cite[Section 6.2]{PecTaq}): in particular, the two integers $k,j$ indicate that, in the definition of $ f  \star_{k}^{j} g$, one identifies $k$ variables in the arguments of $f$ and $g$, among which $j$ are integrated out. In a free setting, we will only need contractions of indices $(k,k)$ and $(k, k-1)$. As in \cite{bianespeicher, knps}, for $1 \leq k \leq m \wedge n$ we shall use the special notation
\begin{eqnarray*}
&& f  \cont{k} g (t_1,\ldots, t_{m+n - 2k}) = f  \star_{k}^{k} g (t_1,\ldots, t_{m+n - 2k})  \\
&& \quad \int_{Z^{k}}f(t_1,\ldots , t_{m-k},s_k, \ldots , s_1)g(s_1, \ldots , s_k , t_{m-k+1},\ldots, t_{m+n - 2k})\mu(ds_1)\cdots \mu(ds_k),
\end{eqnarray*}
and we set moreover $f \cont{0} g =f  \star_{0}^{0} g = f \otimes g$. We observe that, for $k=1,\ldots,m\wedge n$, the star contraction of index $(k,k-1)$ (of $f$ and $g$) is defined by
\begin{eqnarray*}
&& f  \star_{k}^{k-1} g (t_1,\ldots, t_{m+n - 2k+1}) = \\
&& \quad \int_{Z^{k-1}}f(t_1,\ldots , t_{m-k+1}, s_{k-1}, \ldots , s_1) \\
&& \quad\qquad\qquad\qquad\qquad \times g(s_1, \ldots , s_{k-1} , t_{m-k+1},\ldots, t_{m+n - 2k+1})\mu(ds_1)\cdots \mu (ds_{k-1}).
\end{eqnarray*}
\end{itemize}

\bigskip

\noindent In what follows, we shall exploit the useful fact that the class of elementary functions is stable with respect to contractions, in particular: if $f \in \mathscr{E}_m$ and $g \in \mathscr{E}_n$, then $f  \cont{k} g \in \mathscr{E}_{m+n-2k}$ and $f  \star^{k-1}_{k} g \in \mathscr{E}_{m+n-2k+1}$. Also, one can verify the relations
\[
(f  \star_{k}^{k-1} g)^* = g^* \star_{k}^{k-1} f^*,  \quad (f  \cont{k} g)^* = g^* \cont{k} f^*.
\]

\begin{remark}{\rm There exists a crucial difference between `star' and `arc' contractions, that accounts for many technical difficulties one has to deal with in a Poisson framework (both in a classical and free setting), namely: the contraction $\cont{k}$ is a continuous operator from $L^2(Z^m) \times L^2(Z^n) $ into $L^2(Z^{n+m -2k})$ (this can be easily checked by applying the Cauchy--Schwarz inequality), whereas, in general, $\star_{k}^{j}$ (for $k> j$) is not. For instance: one can have that $f_r \to f$ in $L^2(Z^n)$ ($r\to \infty$), without having that $\{ f_r \star_{k}^{k-1} f_r\}$ is a Cauchy sequence in $L^2(Z^{2n-2k+1})$ (actually, $f_r \star_{k}^{k-1} f_r$ need not be square--integrable at all\footnote{Just consider the case $Z=\R$, $n=k=1$, $f_r(x) = x^{-1/4}\mathds{1}_{(0,1/r)}(x)$ and $f =0$.}). For this reason, in the subsequent sections we will be working with kernels that satisfy more regularity assumptions than mere square--integrability (for instance, bounded kernels with bounded support) in order to ensure (e.g., via dominated convergence) some basic continuity properties for star contractions. 
}
\end{remark}

\subsection{Multiple integrals of simple kernels}
We now recall the definition of the multiple stochastic integral of a simple kernel with respect to a centered free measure $\mathfrak{M} = \{\mathfrak{M}(A) : A\in \mathcal{Z}_\mu\}$, defined on an adequate tracial and faithful non--commutative probability space. We assume that 
$$
\mbox{either $\mathfrak{M}=S$ or $\mathfrak{M} = \hat{N}$},
$$
that is: $\mathfrak{M}$ is either a free semicircular measure, or a free centered Poisson measure with Lebesgue control measure ($\mu \equiv {\rm Lebesgue}$). When $Z = \R_+$, and after some straightforward adaptations, our procedure for defining multiple integrals coincides with the construction outlined by Biane and Speicher and Anshelevich, respectively in \cite{bianespeicher} and in \cite{ans0}, in the case of the free Brownian motion and of the free L\'evy processes on the real line. However, the applications we have in mind will lead to new formulae and computations. The reader is referred to the above mentioned papers \cite{ans0, bianespeicher} for more details.

\medskip

\noindent Fix $q\geq 1$, and let $f\in \mathscr{E}_q$ have the form (\ref{e:passa}). The multiple Wigner integral of $f$ with respect to $\hat{N}$ is defined as 
\begin{equation}\label{e:1}
I_q^{\hat{N}} (f_1\otimes \cdots \otimes f_q) = I_q^{\hat{N}} (f) = C_{0, k_1}(\hat{N}(A_1), \mu(A_1))  \cdots C_{0, k_l}(\hat{N}(A_l) , \mu(A_l)),
\end{equation}
where $\{C_{0,k}: k\geq 0\}$ indicates the collection of centered free Charlier polynomials introduced in Section \ref{ss:pre}. Analogously, the multiple Wigner integral of $f$ with respect to $S$ is defined as 
\begin{equation}\label{e:2}
I_q^{S} (f_1\otimes \cdots \otimes f_q) =I_q^{S} (f) = H_{0, k_1}(S(A_1), \mu(A_1))  \cdots H_{0, k_l}(S (A_l) , \mu(A_l)),
\end{equation}
where $\{H_{0,k}: k\geq 0\}$ is the class of generalized Tchebycheff polynomials defined above. We also set $I_0^{\mathfrak{M}}(c) = c\,1$, for every $c\in \mathscr{E}_0 = \mathbb{C}$. Using the recursive properties of Charlier and Tchebycheff polynomials, one sees that our definition of multiple integrals implies the following relation: let $f_1,\ldots,f_q$ be as in formula \eqref{e:passa}, and let $f_0 = \mathds{1}_B$ be such that either $B=A_1$ or $B\cap A_1 = \emptyset$; then
\begin{eqnarray}\label{e:wickmaps}
&& I_{q+1}^{\mathfrak{M}} (f_0\otimes f_1\otimes \cdots \otimes f_q) \\ &&= I_1^{\mathfrak{M}}(f_0)I_q^{\mathfrak{M}}(f_1\otimes \cdots \otimes f_q) 
 - \langle f_0, f_1\rangle_{L^2(Z)}I^{\mathfrak{M}}_{q-1}(f_2\otimes \cdots \otimes f_q ) - \delta_{\mathfrak{M}}\, I^{\mathfrak{M}}_{q}((f_0f_1)\otimes \cdots \otimes f_q ), \nonumber
\end{eqnarray}
where the symbol $\delta_{\mathfrak{M}}$ is defined as
\[
\delta_{{\mathfrak{M}}} = 1 \mbox{ if ${\mathfrak{M}} = \hat{N}$, and } \delta_{{\mathfrak{M}}} = 0 \mbox{ if ${\mathfrak{M}} = S$} .
\]
We extend the definition of $I_q^{\mathfrak{M}}(f)$ to a general $f\in \mathscr{E}_q$ by linearity.

\begin{remark}{\rm 
It is important to notice that relation \eqref{e:wickmaps} implies that the above extension is consistent, in the following sense. If $f$ is an element of $\mathscr{E}_q$ such that $f = \sum_r c_r f_r = \sum_r b_r g_r$, where $c_r,b_r\in \mathbb{C}$ and $\{f_r\}$ and $\{g_r\}$ are two distinct sets of kernels with the form \eqref{e:passa}, then necessarily $I_q^{\mathfrak{M}} (f) =  \sum_r c_r I_q^{\mathfrak{M}}(f_r) = \sum_r b_r I_q^{\mathfrak{M}}(g_r)$. Note also that, by construction, $I_q^{\mathfrak{M}}(f)^* = I_q^{\mathfrak{M}}(f^*)$. }
\end{remark}
A further application of \eqref{e:wickmaps} yields the following fundamental isometric property.
\begin{prop}\label{p:iso} For every $q,q'\geq 0$ and every $f\in \mathscr{E}_q$ and $g\in \mathscr{E}_{q'}$, one has that
\begin{equation}\label{e:iso}
\varphi(I_q^{\mathfrak{M}}(g)^*I^{\mathfrak{M}}_{q'}(f)) = \langle f,g \rangle_{L^2(Z^q)}\,\mathds{1}_{\{q=q'\}}.
\end{equation}
\end{prop}

\noindent Now fix $q\geq 1$. Exploiting the fact that every function in $L^2(Z^q)$ can be approximated in the $L^2$ norm by elements of $\mathscr{E}_q$, we extend the definition of $I_q^{\mathfrak{M}}(f)$ to a general $f\in L^2(Z^q)$, by using the isometric property (\ref{e:iso}). Note that $I_q^{\mathfrak{M}}(f)$ is by construction an element of $L^2(\mathcal{S}({\mathfrak{M}}), \varphi)$, where $\mathcal{S}({\mathfrak{M}})$ is the unital algebra generated by ${\mathfrak{M}}$.

\subsection{Product formulae}\label{ss:prod}

The following statement contains two product formulae for multiple integrals, associated respectively with $\hat{N}$ and $S$. Formula \eqref{e:muls} corresponds to the content of \cite[Proposition 5.3.3]{bianespeicher} (we include a sketch of the proof for the sake of completeness). Formula \eqref{e:mulp} (as well as its extension to general kernels stated in Section \ref{ss:ext}) is new, albeit {\it a posteriori} not surprising -- once it is interpreted in the light of the combinatorial polynomial multiplication formulae proved in \cite{ans1}. See also the subsequent Remark \ref{r:a}.

\begin{theorem}[Product formulae]\label{t:mm} Fix integers $m,n\geq 1$, and let $f \in \mathscr{E}_m$ and $g \in \mathscr{E}_n$. Then,
\begin{eqnarray}
\label{e:mulp}
I^{\hat{N}}_m(f)I^{\hat{N}}_n(g) = \sum_{k=0}^{m \wedge n} I^{\hat{N}}_{m+n-2k}\left( f \cont{k} g\right) +  \sum_{k=1}^{m \wedge n} I^{\hat{N}}_{m+n-2k+1}\left( f \star_{k}^{k-1} g\right), 
\end{eqnarray}
and
\begin{eqnarray}
\label{e:muls}
I^S_m(f)I^S_n(g) = \sum_{k=0}^{m \wedge n} I^S_{m+n-2k}\left( f \cont{k} g\right).
\end{eqnarray}
\end{theorem}

\begin{remark}\label{r:a}{\rm 
\begin{itemize}
\item[(i)] Consider the simple case $f=\mathds{1}^{\otimes m}_A$ and $g=\mathds{1}^{\otimes n}_A$. Then, formulae \eqref{e:mulp}--\eqref{e:muls} together with \eqref{e:1}--\eqref{e:2} yield the following relations:
\begin{eqnarray}
C_{0,m}(\hat{N}(A), \mu(A))C_{0,n}(\hat{N}(A), \mu(A))\!\! &=&\!\! \sum_{k=0}^{n\wedge m} \mu(A)^k C_{0,n+m-2k}(\hat{N}(A), \mu(A)),\label{e:kj}\\
&& + \sum_{k=1}^{n\wedge m} \mu(A)^{k-1} C_{0,n+m-2k+1}(\hat{N}(A), \mu(A))\notag \\
H_{0,m}(S(A), \mu(A))H_{0,n}(S(A), \mu(A)) \!\!&=&\!\! \sum_{k=0}^{n\wedge m} \mu(A)^k H_{0,n+m-2k}(S(A), \mu(A)).\label{e:kjj}
\end{eqnarray}
It is a standard exercise to check that \eqref{e:kj}--\eqref{e:kjj} coincide with the usual combinatorial product formulae for Tchebycheff and free Charlier polynomials -- as proved e.g. in \cite[pp. 122--123]{ans1}.

\item[(ii)] One should compare relation \eqref{e:mulp} with the classical product formulae for Poisson integrals, as proved e.g. in \cite[Section 6.5]{PecTaq}. In particular, if $f,g$ are two symmetric simple kernels as above, and if $I_m^{\hat{\eta}}(f), \, I_n^{\hat{\eta}}(g)$ denote the corresponding multiple Wiener--It\^o integrals with respect to a compensated (classical) Poisson measure with control $\mu$, then one has that
\[
I^{\hat{\eta}}_m(f)I^{\hat{\eta}}_n(g) =  \sum_{k=0}^{m \wedge n}k!\binom{m}{k}\binom{n}{k}\sum_{j=0}^k \binom{k}{j} I^{\hat{\eta}}_{m+n-k-j}\left( \widetilde{f \star_{k}^{j} g}\right), 
\]
where the tilde denotes a symmetrization.
 
\end{itemize}
}
\end{remark}

\medskip

\noindent An immediate consequence of the previous result is the following statement.

\begin{corollary}\label{c:1}{ In both cases ${\mathfrak{M}}=\hat{N}$ or ${\mathfrak{M}}=S$, one has that the class $\{I^{\mathfrak{M}}_q(f) : f\in \mathscr{E}_q\}$ is a unital $\ast$--algebra, with product rule given either by \eqref{e:mulp} or \eqref{e:muls}, and involution $I_q^{\mathfrak{M}}(f)^* = I_q^{\mathfrak{M}}(f^*)$. Moreover, one has that 
\[
\{I^{\mathfrak{M}}_q(f) : f\in \mathscr{E}_q, q\geq 0\} = \mathcal{S}({\mathfrak{M}}),
\]
where $\mathcal{S}({\mathfrak{M}})$ is the unital algebra generated by ${\mathfrak{M}}$, and 
\[
\{I^{\mathfrak{M}}_q(f) : f\in L^2(Z^q), q\geq 0\} = L^2(\mathcal{S}({\mathfrak{M}}), \varphi).
\]
}
\end{corollary}

\subsection{Diagram formulae for simple kernels}\label{ss:dia}

We will now prove new combinatorial formulae for the cumulants associated with free Poisson multiple integrals of simple kernels. To keep the notational complexity of the present paper within reasonable bounds, we choose to state and prove our results only for cumulant of single integrals. Extending the forthcoming results to cumulants of vectors of multiple integrals (as defined e.g. in \cite[p. 175]{NicSpe}) is a standard exercise, that we leave to the interested reader. 

\subsubsection{More on non--crossing partitions}

\noindent We will need some further definitions concerning partitions.

\begin{definition}
\label{respectfulpartitions}{\rm
Let $n_1,\ldots,n_r$ be positive integers such that $n_1 + \cdots + n_r = n$ and consider the partition $\pi^{*} = \{ B_1,\ldots, B_r\}\in NC(n)$, where $B_1 = \left\lbrace 1,\ldots, n_1 \right\rbrace$, $B_2 = \left\lbrace n_1+1,\ldots, n_1 + n_2 \right\rbrace$ and so on until $B_r = \left\lbrace n_1+ \cdots + n_{r-1} + 1,\ldots, n_1 +\cdots + n_r \right\rbrace$. Such a partition $\pi^{*}$ is sometimes called a {\it block partition} and is written $$n_1 \otimes n_2 \otimes \cdots \otimes n_r,$$ see e.g. \cite{knps}. We then say that a partition $\sigma$ of $[n]$ (not necessarily non--crossing) {\it respects} $\pi^{*}$ if no block of $\sigma$ contains more than one element from any given block of $\pi^{*}$. Alternatively, this notion can be formalized by saying that $\sigma$ respects $\pi^{*}$ if and only if $\sigma \wedge \pi^{*} = \hat{0}$. See Fig. 2 for an illustration.}
\end{definition}
\begin{figure}[H]
\begin{center}
\scalebox{1.5} 
{
\begin{pspicture}(0,-0.505)(9.805,0.52)
\definecolor{color1837}{rgb}{0.6,0.8,1}
\definecolor{color1842}{rgb}{0.00392156862745098,0.00392156862745098,0.00392156862745098}
\psdots[dotsize=0.09](3.0,-0.3)
\psdots[dotsize=0.09,linecolor=color1842](3.0,-0.3)
\psdots[dotsize=0.09,linecolor=color1842](3.8,-0.3)
\psdots[dotsize=0.09,linecolor=color1842](1.4,-0.3)
\psdots[dotsize=0.09](0.2,-0.3)
\psdots[dotsize=0.09](0.6,-0.3)
\psdots[dotsize=0.09](1.0,-0.3)
\psdots[dotsize=0.09](2.0,-0.3)
\psdots[dotsize=0.09](2.4,-0.3)
\psdots[dotsize=0.09](3.4,-0.3)
\psdots[dotsize=0.09](4.4,-0.3)
\psdots[dotsize=0.09](0.2,-0.3)
\psdots[dotsize=0.09](0.6,-0.3)
\psdots[dotsize=0.09](3.8,-0.3)
\psdots[dotsize=0.09](4.4,-0.3)
\psdots[dotsize=0.09](1.0,-0.3)
\psdots[dotsize=0.09](1.4,-0.3)
\psdots[dotsize=0.09](2.0,-0.3)
\psdots[dotsize=0.09](2.4,-0.3)
\psdots[dotsize=0.09](3.0,-0.3)
\psdots[dotsize=0.09](3.4,-0.3)
\psline[linewidth=0.04cm,linecolor=color1837](2.4,0.1)(2.4,-0.3)
\psline[linewidth=0.04cm](2.0,-0.1)(2.0,-0.3)
\psdots[dotsize=0.09](8.2,-0.3)
\psdots[dotsize=0.09,linecolor=color1842](8.2,-0.3)
\psdots[dotsize=0.09,linecolor=color1842](9.0,-0.3)
\psdots[dotsize=0.09,linecolor=color1842](6.6,-0.3)
\psdots[dotsize=0.09](5.4,-0.3)
\psdots[dotsize=0.09](5.8,-0.3)
\psdots[dotsize=0.09](6.2,-0.3)
\psdots[dotsize=0.09](7.2,-0.3)
\psdots[dotsize=0.09](7.6,-0.3)
\psdots[dotsize=0.09](8.6,-0.3)
\psdots[dotsize=0.09](9.6,-0.3)
\psdots[dotsize=0.09](5.4,-0.3)
\psdots[dotsize=0.09](5.8,-0.3)
\psdots[dotsize=0.09](9.0,-0.3)
\psdots[dotsize=0.09](9.6,-0.3)
\psdots[dotsize=0.09](6.2,-0.3)
\psdots[dotsize=0.09](6.6,-0.3)
\psdots[dotsize=0.09](7.2,-0.3)
\psdots[dotsize=0.09](7.6,-0.3)
\psdots[dotsize=0.09](8.2,-0.3)
\psdots[dotsize=0.09](8.6,-0.3)
\psline[linewidth=0.04cm](3.0,-0.1)(1.4,-0.1)
\psline[linewidth=0.04cm](1.4,-0.1)(1.4,-0.3)
\psline[linewidth=0.04cm](3.0,-0.1)(3.0,-0.3)
\psline[linewidth=0.04cm](0.6,-0.3)(0.6,0.3)
\psline[linewidth=0.04cm](0.6,0.3)(3.8,0.3)
\psline[linewidth=0.04cm](3.8,0.3)(3.8,-0.3)
\psline[linewidth=0.04cm](4.4,-0.3)(4.4,0.5)
\psline[linewidth=0.04cm](4.4,0.5)(0.2,0.5)
\psline[linewidth=0.04cm](0.2,0.5)(0.2,-0.3)
\psline[linewidth=0.01cm](0.0,-0.5)(1.6,-0.5)
\psline[linewidth=0.01cm](1.6,-0.3)(1.6,-0.5)
\psline[linewidth=0.01cm](0.0,-0.3)(0.0,-0.5)
\psline[linewidth=0.01cm](1.8,-0.3)(1.8,-0.5)
\psline[linewidth=0.01cm](1.8,-0.5)(2.6,-0.5)
\psline[linewidth=0.01cm](2.6,-0.3)(2.6,-0.5)
\psline[linewidth=0.01cm](2.8,-0.3)(2.8,-0.5)
\psline[linewidth=0.01cm](2.8,-0.5)(4.0,-0.5)
\psline[linewidth=0.01cm](4.0,-0.3)(4.0,-0.5)
\psline[linewidth=0.01cm](4.2,-0.3)(4.2,-0.5)
\psline[linewidth=0.01cm](4.2,-0.5)(4.6,-0.5)
\psline[linewidth=0.01cm](4.6,-0.5)(4.6,-0.3)
\psdots[dotsize=0.09](8.2,-0.3)
\psdots[dotsize=0.09,linecolor=color1842](8.2,-0.3)
\psdots[dotsize=0.09,linecolor=color1842](9.0,-0.3)
\psdots[dotsize=0.09,linecolor=color1842](6.6,-0.3)
\psdots[dotsize=0.09](5.4,-0.3)
\psdots[dotsize=0.09](5.8,-0.3)
\psdots[dotsize=0.09](6.2,-0.3)
\psdots[dotsize=0.09](7.2,-0.3)
\psdots[dotsize=0.09](7.6,-0.3)
\psdots[dotsize=0.09](8.6,-0.3)
\psdots[dotsize=0.09](9.6,-0.3)
\psdots[dotsize=0.09](5.4,-0.3)
\psdots[dotsize=0.09](5.8,-0.3)
\psdots[dotsize=0.09](9.0,-0.3)
\psdots[dotsize=0.09](9.6,-0.3)
\psdots[dotsize=0.09](6.2,-0.3)
\psdots[dotsize=0.09](6.6,-0.3)
\psdots[dotsize=0.09](7.2,-0.3)
\psdots[dotsize=0.09](7.6,-0.3)
\psdots[dotsize=0.09](8.2,-0.3)
\psdots[dotsize=0.09](8.6,-0.3)
\psline[linewidth=0.01cm](5.2,-0.5)(6.8,-0.5)
\psline[linewidth=0.01cm](6.8,-0.3)(6.8,-0.5)
\psline[linewidth=0.01cm](5.2,-0.3)(5.2,-0.5)
\psline[linewidth=0.01cm](7.0,-0.3)(7.0,-0.5)
\psline[linewidth=0.01cm](7.0,-0.5)(7.8,-0.5)
\psline[linewidth=0.01cm](7.8,-0.3)(7.8,-0.5)
\psline[linewidth=0.01cm](8.0,-0.3)(8.0,-0.5)
\psline[linewidth=0.01cm](8.0,-0.5)(9.2,-0.5)
\psline[linewidth=0.01cm](9.2,-0.3)(9.2,-0.5)
\psline[linewidth=0.01cm](9.4,-0.3)(9.4,-0.5)
\psline[linewidth=0.01cm](9.4,-0.5)(9.8,-0.5)
\psline[linewidth=0.01cm](9.8,-0.5)(9.8,-0.3)
\psline[linewidth=0.04cm](9.6,0.5)(5.4,0.5)
\psline[linewidth=0.04cm](7.2,-0.1)(6.6,-0.1)
\psline[linewidth=0.04cm](6.6,-0.1)(6.6,-0.3)
\psline[linewidth=0.04cm](7.2,-0.1)(7.2,-0.3)
\psline[linewidth=0.04cm](8.2,-0.1)(7.6,-0.1)
\psline[linewidth=0.04cm](7.6,-0.1)(7.6,-0.3)
\psline[linewidth=0.04cm](8.2,-0.1)(8.2,-0.3)
\psline[linewidth=0.04cm](8.6,0.1)(6.2,0.1)
\psline[linewidth=0.04cm](8.6,0.1)(8.6,-0.3)
\psline[linewidth=0.04cm](6.2,0.1)(6.2,-0.3)
\psline[linewidth=0.04cm](9.0,-0.3)(9.0,0.3)
\psline[linewidth=0.04cm](9.0,0.3)(5.8,0.3)
\psline[linewidth=0.04cm](5.8,0.3)(5.8,-0.3)
\psline[linewidth=0.04cm](5.4,0.5)(5.4,-0.3)
\psline[linewidth=0.04cm](9.6,0.5)(9.6,-0.3)
\psline[linewidth=0.04cm,linecolor=color1837](1.0,0.1)(3.4,0.1)
\psline[linewidth=0.04cm,linecolor=color1837](1.0,0.1)(1.0,-0.3)
\psline[linewidth=0.04cm,linecolor=color1837](3.4,0.1)(3.4,-0.3)
\psdots[dotsize=0.09](1.0,-0.3)
\psdots[dotsize=0.09](3.4,-0.3)
\psdots[dotsize=0.09](2.4,-0.3)
\end{pspicture} 
}
\caption{\sl A crossing (left) and a non--crossing (right) partition of $[10]$ that respect the block partition $4\otimes 2 \otimes 3 \otimes 1$.}
\end{center}
\end{figure}
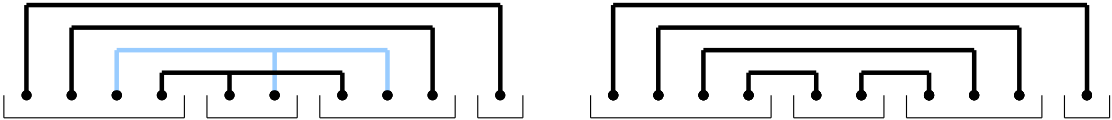
\begin{remark}{\rm
\label{remarkdefalternnoncrossing} Let the notation and assumptions of Lemma \ref{l:1} prevail, and assume that $\pi\in NC(n)$ is respectful of the block partition $\pi^{*}$, and has blocks of size at least 2. Then, one can choose the interval block $V$ to have size exactly $2$, more precisely: $\pi$ contains necessarily a block of the type $V = \left\lbrace k,k+1 \right\rbrace $ for some $1 \leq k \leq n$, and the integers $k$ and $k+1$ belong to two distinct consecutive blocks of $\pi^*$.}
\end{remark}

\noindent To conclude, we now define the notion of partition $\sigma$ connecting a block partition $\pi^{*}$.

\begin{definition}{\rm 
\label{connectedpartitions}
Let $\sigma \in \mathcal{P}(n)$ and let $\pi^{*} = n_1 \otimes n_2 \otimes \cdots \otimes n_r$, where, as before, $n_1,\ldots, n_r$ are positive integers such that $n_1+\cdots + n_r = n$. Two blocks $B_1, B_2$ in $\pi^{*}$ are said to be \textit{linked} by $\sigma$ if there is a block in $\sigma$ containing at least an element of $B_1$ and at least an element of $B_2$. Denote by $C_{\sigma}$ the graph whose vertices are the blocks of $\pi{*}$ and that has an edge between $B_1$ and $B_2$ if and only if $\sigma$ links $B_1$ to $B_2$. The partition $\sigma$ is said to {\it connect} $\pi^{*}$ if the graph $C_{\sigma}$ is connected. Alternatively, this notion can be formalized by saying that $\sigma \in \mathcal{P}(n)$ connects $\pi^{*}$ if and only if $\sigma \vee \pi^{*} = \hat{1}$.}
\end{definition}
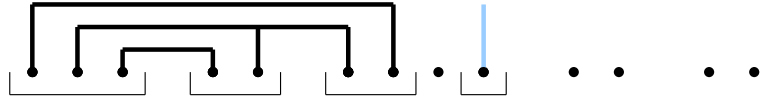
\begin{figure}[H]
\begin{center}
\scalebox{1.5} 
{
\begin{pspicture}(0,-0.405)(4.405,0.42)
\definecolor{color1842}{rgb}{0.00392156862745098,0.00392156862745098,0.00392156862745098}
\definecolor{color2411}{rgb}{0.6,0.8,1.0}
\psdots[dotsize=0.09](3.0,-0.2)
\psdots[dotsize=0.09,linecolor=color1842](3.0,-0.2)
\psdots[dotsize=0.09](0.2,-0.2)
\psdots[dotsize=0.09](0.6,-0.2)
\psdots[dotsize=0.09](1.0,-0.2)
\psdots[dotsize=0.09](1.8,-0.2)
\psdots[dotsize=0.09](3.4,-0.2)
\psdots[dotsize=0.09](0.2,-0.2)
\psdots[dotsize=0.09](0.6,-0.2)
\psdots[dotsize=0.09](1.0,-0.2)
\psdots[dotsize=0.09](1.8,-0.2)
\psdots[dotsize=0.09](2.2,-0.2)
\psdots[dotsize=0.09](3.0,-0.2)
\psdots[dotsize=0.09](3.4,-0.2)
\psdots[dotsize=0.09](3.0,-0.2)
\psdots[dotsize=0.09,linecolor=color1842](3.0,-0.2)
\psdots[dotsize=0.09](0.2,-0.2)
\psdots[dotsize=0.09](0.6,-0.2)
\psdots[dotsize=0.09](1.0,-0.2)
\psdots[dotsize=0.09](1.8,-0.2)
\psdots[dotsize=0.09](2.2,-0.2)
\psdots[dotsize=0.09](3.4,-0.2)
\psdots[dotsize=0.09](0.2,-0.2)
\psdots[dotsize=0.09](0.6,-0.2)
\psdots[dotsize=0.09](4.2,-0.2)
\psdots[dotsize=0.09](1.0,-0.2)
\psdots[dotsize=0.09](1.8,-0.2)
\psdots[dotsize=0.09](2.2,-0.2)
\psdots[dotsize=0.09](3.0,-0.2)
\psdots[dotsize=0.09](3.4,-0.2)
\psline[linewidth=0.01cm](0.0,-0.4)(1.2,-0.4)
\psline[linewidth=0.01cm](1.2,-0.2)(1.2,-0.4)
\psline[linewidth=0.01cm](0.0,-0.2)(0.0,-0.4)
\psline[linewidth=0.01cm](1.6,-0.2)(1.6,-0.4)
\psline[linewidth=0.01cm](1.6,-0.4)(2.4,-0.4)
\psline[linewidth=0.01cm](2.4,-0.2)(2.4,-0.4)
\psline[linewidth=0.01cm](2.8,-0.2)(2.8,-0.4)
\psline[linewidth=0.01cm](2.8,-0.4)(3.6,-0.4)
\psline[linewidth=0.01cm](3.6,-0.2)(3.6,-0.4)
\psline[linewidth=0.01cm](4.0,-0.2)(4.0,-0.4)
\psline[linewidth=0.01cm](4.0,-0.4)(4.4,-0.4)
\psline[linewidth=0.01cm](4.4,-0.4)(4.4,-0.2)
\psline[linewidth=0.04cm](1.8,0.0)(1.0,0.0)
\psline[linewidth=0.04cm](1.8,0.0)(1.8,-0.2)
\psline[linewidth=0.04cm](2.2,0.2)(2.2,-0.2)
\psline[linewidth=0.04cm](3.0,0.2)(0.6,0.2)
\psline[linewidth=0.04cm](3.0,0.2)(3.0,-0.2)
\psline[linewidth=0.04cm](0.6,0.2)(0.6,-0.2)
\psline[linewidth=0.04cm](0.2,0.4)(0.2,-0.2)
\psline[linewidth=0.04cm,linecolor=color2411](4.2,0.4)(4.2,-0.2)
\psline[linewidth=0.04cm](1.0,-0.2)(1.0,0.0)
\psline[linewidth=0.04cm](3.4,-0.2)(3.4,0.4)
\psline[linewidth=0.04cm](0.2,0.4)(3.4,0.4)
\psdots[dotsize=0.09](4.2,-0.2)
\end{pspicture} 
}
\caption{\sl A partition of ${NC}(8)$ that respects but does not connect $3\otimes 2 \otimes 2 \otimes 1$.}
\end{center}
\end{figure}

\noindent Fix two positive integers $q$ and $m$ such that $m,q \geq 1$. From now on, and unless otherwise specified, $\pi^{*}$ will denote the block partition of $[mq]$ given by $$\pi^{*} := \underbrace{q \otimes\cdots \otimes q}_{m \mbox{\tiny{ times}}}.$$

\begin{definition}{\rm
\label{definitionsetsNC}
Let ${NC}\left(\left[mq \right]  , \pi^{*} \right) $ and ${NC}^{0}\left(\left[mq \right]  , \pi^{*} \right) $ be the sets of partitions of $\left[mq \right] $ defined by 
\begin{eqnarray*}
&& {NC}\left(\left[mq \right] , \pi^{*} \right) = \left\lbrace \sigma \in {NC}(mq) : \sigma \vee \pi^{*} = \hat{1} \rm{\ and\ } \sigma \wedge \pi^{*} = \hat{0} \right\rbrace\\
&& {NC}^{0}\left(\left[mq \right] , \pi^{*} \right) = \left\lbrace \sigma \in {NC}(mq) : \sigma \wedge \pi^{*} = \hat{0} \right\rbrace  .
\end{eqnarray*}
Furthermore, we define the four following subsets of ${NC}\left(\left[mq \right] , \pi^{*} \right) $ and ${NC}^{0}\left(\left[mq \right] , \pi^{*} \right) $: 
\begin{eqnarray*}
&& {NC}_{2}\left(\left[mq \right] , \pi^{*} \right) = \left\lbrace \sigma \in {NC}\left(\left[mq \right] , \pi^{*} \right) : \vert b \vert = 2, \forall b \in \sigma \right\rbrace \\
&& {NC}_{\geq 2}\left(\left[mq \right] , \pi^{*} \right) = \left\lbrace \sigma \in {NC}\left(\left[mq \right] , \pi^{*} \right) : \vert b \vert \geq 2, \forall b \in \sigma \right\rbrace\\ 
&& {NC}^{0}_{2}\left(\left[mq \right] , \pi^{*} \right) = \left\lbrace \sigma \in {NC}^{0}\left(\left[mq \right] , \pi^{*} \right) : \vert b \vert = 2, \forall b \in \sigma \right\rbrace \\
&& {NC}^{0}_{\geq 2}\left(\left[mq \right] , \pi^{*} \right) = \left\lbrace \sigma \in {NC}^{0}\left(\left[mq \right] , \pi^{*} \right) : \vert b \vert \geq 2, \forall b \in \sigma \right\rbrace.
\end{eqnarray*}
}
\end{definition}    

\begin{remark}\label{r:riordan}{\rm When specialized to the case $q=1$, some of the previous classes of partitions reduce to simple objects. In particular, one has that ${NC}\left(\left[m\right]  , \pi^{*} \right) = \{\hat{1}\}$, ${NC}^{0}\left(\left[m \right]  , \pi^{*} \right) = NC(m)$ and ${NC}^{0}_{\geq 2}\left(\left[m \right] , \pi^{*} \right)$ equals the collection of all non--crossing partitions of $[m]$ having no singletons. This implies that $\left|{NC}^{0}_{\geq 2}\left(\left[m \right] , \pi^{*} \right)\right| = R_m$, where $R_m$ denotes the $m^{\mbox{\tiny{th}}}$ {\it Riordan number} (see \cite{ber, np}). Observe also that ${NC}_{2}\left(\left[m \right] , \pi^{*} \right) = \emptyset$, for $m\geq 3$.
}
\end{remark}

\subsubsection{Statements}

\noindent The next statement contains the announced diagram formulae for multiple integrals. Note that relations \eqref{e:dp}--\eqref{e:srp} are new, whereas \eqref{e:ds}--\eqref{e:ds2} are basically equivalent to \cite[Proposition 1.38]{knps} (a sketch of the proof is once again included for the sake of completeness). The Haagerup--type inequality \eqref{e:srs} is proved in \cite[Theorem 5.3.4]{bianespeicher}. In the proof, we shall use the definition of the class $\mathscr{E}_q^{00}$, of purely non--diagonal elements, as introduced in Lemma \ref{l:taborn}.

\begin{definition}[Partitions and tensors]\label{d:tensor}{\rm Let $q,m\geq 1$, and consider a function $f$ in $q$ variables. Given a partition $\sigma$ of $[mq]$, we define the function $f_{\sigma}$, in $\vert\sigma\vert$ variables, as the mapping obtained by identifying the variables $x_i$ and $x_j$ in the argument of the tensor
\begin{equation}
\label{formegeneraledefsigma}
f\otimes \ldots \otimes f\left(x_1,\ldots,x_{mq} \right)  = \prod_{j=1}^{m}f\left(x_{(j-1)q + 1},\ldots,x_{jq} \right)
\end{equation}
if and only if $i$ and $j$ are in the same block of $\sigma$. For instance, if $q=2,\, m=4$ and $\sigma = \{\{1,8\}, \{2,3\}, \{4,5\}, \{6,7\}\}\in NC_2([8], \pi^*)$, then
\[
f_{\sigma} (x,y,v,w) = f(x,y)f(y,v)f(v,w)f(w,x).
\]
}
\end{definition}

\begin{theorem}[Diagram formulae \& spectral bounds]
\label{freediagramformulae}
Let ${\mathfrak{M}}$ be a centered free random measure on $(Z, \mathcal{Z})$, with non--atomic control $\mu$. For any $f \in \mathscr{E}_q$ and any integer $q\geq 1$ and $m\geq 2$, it holds that:
\begin{itemize}
\item[\rm (i)] If ${\mathfrak{M}} = \hat{N}$ is a centered free Poisson random measure, then
\begin{eqnarray}\label{e:dp}
\kappa_{m}\left( I_{q}^{\hat{N}}(f)\right)  &=& \sum_{\sigma \in {NC}_{\geq 2}\left(\left[mq \right], \pi^{*}  \right) }\int_{Z^{\vert \sigma \vert}}f_{\sigma}d\mu^{\vert \sigma \vert}\\
\varphi\left( I_{q}^{\hat{N}}(f)^m\right)  &=& \sum_{\sigma \in {NC}^{0}_{\geq 2}\left(\left[mq \right], \pi^{*}  \right) }\int_{Z^{\vert \sigma \vert}}f_{\sigma}d\mu^{\vert \sigma \vert}.\label{e:dp2}
\end{eqnarray}
Moreover, let $I_{q}^{\hat{N}}(f)$ is self--adjoint, and let $B\subset Z$ and $D\in (0,\infty)$ be, respectively, (a) a measurable set such that $\mu(B) := K<\infty$ and the support of $f$ is contained in $B\times \cdots \times B$, and (b) a constant satisfying $|f| \leq D$: then, the spectral radius of $I_{q}^{\hat{N}}(f)$ satisfies the following inequality
\begin{equation}\label{e:srp}
\rho\left(I_{q}^{\hat{N}}(f)\right) \leq 4^q\,  \max\{1; DK\}^{q/2}.
\end{equation}
\item[\rm (ii) ] If ${\mathfrak{M}} = S$ is a centered semicircular measure, then
\begin{eqnarray}\label{e:ds}
\kappa_{m}\left( I_{q}^{S}(f)\right)  &=& \sum_{\sigma \in {NC}_{2}\left(\left[mq \right], \pi^{*}  \right) }\int_{Z^{mq/2}}f_{\sigma}d\mu^{mq/2}\\
\label{e:ds2} \varphi\left( I_{q}^{S}(f)^m\right)  &=& \sum_{\sigma \in {NC}^{0}_{2}\left(\left[mq \right], \pi^{*}  \right) }\int_{Z^{mq/2}}f_{\sigma}d\mu^{mq/2}.
\end{eqnarray}
Moreover, if $f$ is mirror--symmetric, one has the spectral bound
\begin{equation}\label{e:srs}
\rho\left(I_{q}^{S}(f)\right) \leq (q+1)\| f\|_{L^2(Z^q)}.
\end{equation}

\end{itemize}
\end{theorem}

\begin{remark}
{\rm Albeit largely sufficient for our needs, the spectral bound \eqref{e:srp} is not satisfying, in particular because the quantity $\max\{1;DK\}$ is not converging to zero whenever $DK\to 0$. By inspection of the forthcoming proof, one sees that, in order to obtain a spectral bound verifying such a basic continuity property, one would need to show that the minimal number of blocks of a partition in $NC_{\geq 2}^0([2mq],\pi^*)$ strictly increases as an affine function of $m$ (in the proof we use the trivial bound: $|\sigma|\geq q$ for every $\sigma\in NC_{\geq 2}^0([2mq],\pi^*)$). We prefer to think of this issue as a separate problem, and we leave it open for further research.}
\end{remark}

\subsection{Extension to general kernels}\label{ss:ext}

We now present two statements, showing that one can extend the results proved above to more general kernels. The proofs --- that are standard and left to the reader --- follow from Lemma \ref{l:lot} and Corollary \ref{c:1} and, respectively, from Lemma \ref{l:taborn} and from the fact that $\mathscr{E}_q^{00}$ is dense in $L^2(Z^q)$. 

\begin{prop}\label{p:taxi1}
\begin{itemize}
\item[\rm (i)] Denote by $L_b(\mathcal{S}(\hat{N}),\varphi)$ the collection of all objects of the type $I_q^{\hat{N}}(f)$, where $f$ is a bounded function with bounded support (so that $L_b(\mathcal{S}(\hat{N}),\varphi)\, \subset L^2(\mathcal{S}(\hat{N}),\varphi)$). Then, $L_b(\mathcal{S}(\hat{N}),\varphi)$ is a unital $\ast$--algebra, with involution $I_q^{\hat{N}}(f)^* = I_q^{\hat{N}}(f^*)$, and product rule given by formula \eqref{e:mulp}. The trace state $\varphi$ on the class \[\mathcal{S}(\hat{N}) =\{I^{\hat{N}}_q(f) : f\in \mathscr{E}_q, q\geq 0\}\] extends to $L_b(\mathcal{S}(\hat{N}),\varphi)$, in such a way that the cumulant and moment formulae \eqref{e:dp}--\eqref{e:dp2} and the spectral bound \eqref{e:srp} continue to hold.

\item[\rm (ii)] The collection of all objects of the type $I_q^{S}(f)$, where $f\in L^2(Z^q)$, is a unital $\ast$--algebra coinciding with $L^2(\mathcal{S}(S), \varphi)$, with involution $I_q^{S}(f)^* = I_q^{S}(f^*)$, and product rule given by formula \eqref{e:muls}. The trace state $\varphi$ on the class $\mathcal{S}(S)= \{I^S_q(f) : f\in \mathscr{E}_q, q\geq 0\}$ extends to $L^2(\mathcal{S}(S), \varphi)$, in such a way that the cumulant and moment formulae \eqref{e:ds}--\eqref{e:ds2} and the spectral bound \eqref{e:srs} continue to hold.
\end{itemize}
\end{prop}

\begin{remark}[On single integrals]{\rm Let $I^{\hat{N}}_1(f) \in L_b(\mathcal{S}(\hat{N}),\varphi)$ be such that $f=f^*$. Then, Proposition \ref{p:taxi1}--(i) together with Remark \ref{r:riordan} imply the following relations: for every integer $m\geq 2$,
\begin{eqnarray}
\kappa_m(I^{\hat{N}}_1(f)) &=& \int_Z f^m d\mu, \label{e:1a}\\
\varphi(I^{\hat{N}}_1(f)^m) &=& \sum_{\sigma} \prod_{b\in \sigma} \int_Z f^{|b|}d\mu,\label{e:1b}
\end{eqnarray}
where the sum runs over the class of all non--crossing partitions of $[m]$ having blocks of size at least 2. In particular, if $f = \mathds{1}_A$, where $\mu(A) = \lambda\in (0,\infty)$, then $I^{\hat{N}}_1(f) = \hat{N}(A)\sim \hat{P}(\lambda)$, and one recovers from \eqref{e:1b} the content of \cite[Proposition 2.4]{np}, according to which
\[
\varphi(\hat{N}(A)^m) = \sum_{j=1}^m \lambda^j R_{m,j},
\]
where $R_{m,j}$ denotes the number of non--crossing partitions of $[m]$ having exactly $j$ blocks, and such that each block has at least size $2$.}
\end{remark}

\noindent The next statement, which is a direct consequence of \eqref{e:1a}, contains some simple facts about the semicircular approximation of random variables of the type $I^{\hat{N}}_1(f)$ (the straightforward proof is left to the reader). In the next section, we will provide and exhaustive generalization to free multiple integrals of arbitrary order. 

\begin{prop} \label{p:noc1}
\begin{itemize}

\item[\rm (i)] Let $I^{\hat{N}}_1(f) \in L_b(\mathcal{S}(\hat{N}),\varphi)$ be such that $f\neq 0$ and $f = f^*$. Then, $I^{\hat{N}}_1(f)$ cannot have a semicircular distribution.

\item[\rm (ii)] Let $\{ I^{\hat{N}}_1(f_n) : n \geq 1\} \subset L_b(\mathcal{S}(\hat{N}),\varphi)$ be such that $f_n = f_n^*$, $\|f_n\|^2_{L^2(Z)} \to \alpha^2$ and $\sup_n \int_Z |f_n|^p d\mu <\infty$ for every $p\geq 2$. Then $I^{\hat{N}}_1(f_n)$ converges in law to $\mathcal{S}(0,\alpha^2)$ if and only if $\int_Z f_n^4 d\mu \to 0$. 

\end{itemize}
\end{prop}

\section{Semicircular limits for free Poisson multiple integrals}

In this section, we prove semicircular limit theorems for free Poisson multiple integrals. In order to obtain neater statements, we will focus on sequences that are \textit{tamed}, in a sense to be specified in the next definition.  Recall that we work under the convention that $Z = \mathbb{R}^d$ and $\mu$ equals the Lebesgue measure.

\begin{definition}[Tamed sequences]\label{deftamedsuites}{\rm Let $q\geq 1$. We say that the sequence $\{g_n : n\geq 1\}\subset L^2(Z^q)$ is {\it tamed} if the following conditions hold: every $g_n$ is bounded and has bounded support and, for every $m\geq 2$ and every $\sigma\in \mathcal{P}(mq)$ such that $\sigma\wedge \pi^* = \hat{0}$, the numerical sequence 
\begin{equation}
\label{tamedsequencebornee}
\int_{Z^{|\sigma|}} \left| g_n \right|_\sigma d\mu^{|\sigma|}, \quad n\geq 1,
\end{equation}
is bounded, where $\pi^*\in \mathcal{P}(mq)$ is the block partition with $m$ consecutive blocks of size $q$, and the function $\left| g_n \right|_\sigma$, in $|\sigma|$ variables, is defined according to Definition \ref{d:tensor} in the case $f = |g_n|$. Note that the condition $\sigma\wedge \pi^* = \hat{0}$ implies that, necessarily, $|\sigma|\geq q$.

}\label{d:tamed}
\end{definition}

\noindent The next statement provides useful sufficient conditions in order for a sequence $\{f_n\}$ to be tamed: this basically consists in requiring that $\{f_n\}$ concentrates asymptotically, without exploding, around a hyperdiagonal.

\begin{lemma} \label{l:tamed}Fix $q\geq 2$, and consider a sequence $\{f_n : n\geq 1\}\subset L^2(Z^q)$. Assume that there exist strictly positive numerical sequences $\{M_n, z_n, \alpha_n: n\geq 1\}$ such that $\alpha_n/z_n \to 0$ and the following properties are satisfied:
\begin{itemize}
\item[\rm (a)] The support of $f_n$ is contained in the set $(-z_n,z_n)^d\times\cdots\times (-z_n, z_n)^d$ (Cartesian product of order $q$);

\item[\rm (b)] $|f_n|\leq M_n$; 

\item[\rm (c)] $f_n(x_1,\ldots,x_q) = 0$, whenever there exist $x_i, x_j$ such that $\| x_i - x_j\|_{\mathbb{R}^d}>\alpha_n$.
\item[\rm (d)] For every integer $m\geq q$, the mapping $n\mapsto M_n^m z_n^{d} (\alpha_n^d)^{m-1}$ is bounded.
\end{itemize}
Then, $\left\lbrace f_n \colon n\geq 1 \right\rbrace $ is tamed.
\end{lemma}
\noindent{\bf Proof}. Fix $m\geq 2$, as well as $\sigma\in \mathcal{P}(mq)$ such that $\sigma\wedge \pi^* = \hat{0}$. By definition,
\begin{eqnarray*}
&&\int_{Z^{|\sigma|}} |f_n|_\sigma d\mu^{|\sigma|}\leq M_n^{|\sigma|} \int_{(-z_n,z_n)^d}\cdots \int_{(-z_n,z_n)^d} \mathds{1}_{\{\|x_i - x_j\|_{\R^d}\leq|\sigma|\alpha_n\} }\mu(dx_1)\cdots \mu(dx_{|\sigma|})\\
&& = M_n^{|\sigma|}z_n^{|\sigma|d} \int_{(-1,1)^d}\cdots \int_{(-1,1)^d} \mathds{1}_{\{\|y_i - y_j\|_{\R^d}\leq|\sigma|\alpha_n/z_n\} }\mu(dy_1)\cdots \mu(dy_{|\sigma|}).
\end{eqnarray*}
Write $|\sigma|\alpha_n/z_n = r_n$. Applying the change of variables $v_1 = x_1$, $y_i = v_1+ r_n v_i$, for $i=2,\ldots,q$, yields that the previous expression is asymptotically equivalent to
\[
M_n^{|\sigma|}z_n^{|\sigma|d}(r_n^d)^{|\sigma|-1} \times \int_{\R^d}\cdots \int_{\R^d}\mathds{1}_{\{\| v_j\|_{\R^d}\leq1, \,\, \forall j \geq 1\} } \mathds{1}_{\{\|v_i - v_j\|_{\R^d}\leq1, \,\, \forall i,j\geq 2\} }\mu(dv_1)\cdots \mu(dv_{|\sigma|}),
\]
from which we infer the desired conclusion.
\qed

\medskip

\noindent The next result is one of the main achievements of the paper.

\begin{theorem}
\label{fourthcumulant}
Fix $q\geq 2$, and consider a tamed sequence $\left\lbrace f_{n} : n\geq 1 \right\rbrace$ of mirror symmetric functions such that $\underset{n \rightarrow \infty}{\rm lim} \Vert f_n \Vert_{L^{2}\left( Z^q\right) }^{2} =\alpha^2 < \infty$. Then, the following four conditions are equivalent, as $n\to \infty$
\begin{itemize}

\item[\rm (i)] $I_{q}^{\hat{N}}\left( f_{n}\right)$ converges in law to $\mathcal{S}(0, \alpha^2)$;

\item[\rm (ii)] $\kappa_4(I_{q}^{\hat{N}}\left( f_{n}\right))\to 0$;

\item[\rm (iii)] $\varphi\left(\left[ I_{q}^{\hat{N}}\left( f_{n}\right) \right] ^4 \right) {\to} 2\alpha^4$;

\item[\rm (iv)] $\| f_n \cont{k} f_n\|_{L^2(Z^{2q-2k})} \rightarrow 0$ for all $k \in \left\lbrace 1,\ldots,q-1 \right\rbrace $ and $\left\| f_n \star_{k}^{k-1} f_n \right\|_{L^2(Z^{2q-2k+1})} \rightarrow 0$ for all $k \in \left\lbrace 2,\ldots,q \right\rbrace$.
\end{itemize}

\end{theorem}  

\begin{remark}{\rm 
\begin{itemize}

\item[(a)] Recall from Section \ref{ss:pre} that, according to our terminology, stating that $I_{q}^{\hat{N}}\left( f_{n}\right)$ converges in law to $\mathcal{S}(0, \alpha^2)$ is the same as requiring that the moments of $I_{q}^{\hat{N}}\left( f_{n}\right)$ all converge to the corresponding moments of $\mathcal{S}(0, \alpha^2)$: it follows that the implications (i) $\Longrightarrow$ (ii), (iii) in the previous statement are just a direct consequence of our definition. We stress that the statement would not hold, in general, if (i) was replaced by the condition that the spectral measure of $I_{q}^{\hat{N}}\left( f_{n}\right)$ weakly converges to that of $\mathcal{S}(0,\alpha^2)$ (in this case, one would need e.g. some additional uniform control on the spectral radius of $I_{q}^{\hat{N}}\left( f_{n}\right)$).
\item[(b)] Observe the following basic identity, valid for every mirror--symmetric $f\in L^2(Z^q)$ that is bounded and has bounded support (the proof is based on a standard application of Fubini theorem):
\[
\|f\star_1^0 f\|_{L^2(Z^{2q-1})} = \|f\star_k^{k-1} f\|_{L^2(Z)}. 
\]

\end{itemize}
}
\end{remark}

\noindent{\bf Proof of Theorem \ref{fourthcumulant}}. The implications (i) $\Longrightarrow$ (ii) $\Longleftrightarrow$ (iii) are trivial. In order to prove the implication (iii) $\Longrightarrow$ (iv), observe that, according to \eqref{e:mulp}, 
\begin{eqnarray*}
I^{\hat{N}}_q(f_n)^2 &=& \sum_{k=0}^{q} I^{\hat{N}}_{2q-2k}\left( f_n \cont{k} f_n\right) +  \sum_{k=1}^{q} I^{\hat{N}}_{2q-2k+1}\left( f_n \star_{k}^{k-1} f_n\right) \\
&=& f_n \cont{q} f_n + I^{\hat{N}}_{2q}(f_n \otimes f_n) + \sum_{k=1}^{q-1} I^{\hat{N}}_{2q-2k}\left( f_n \cont{k} f_n\right) +  \sum_{k=1}^{q} I^{\hat{N}}_{2q-2k+1}\left( f_n \star_{k}^{k-1} f_n\right) \\
&=& \Vert f_n \Vert_{L^{2}\left( Z^q\right) }^{2} + I^{\hat{N}}_{2q}(f_n \otimes f_n) + \sum_{k=1}^{q-1} I^{\hat{N}}_{2q-2k}\left( f_n \cont{k} f_n\right) +  \sum_{k=1}^{q} I^{\hat{N}}_{2q-2k+1}\left( f_n \star_{k}^{k-1} f_n\right).
\end{eqnarray*}
Using the isometric properties of multiple integrals \eqref{e:iso}, together with the fact that multiple integrals of different orders are orthogonal in $L^2(\mathscr{A}, \varphi)$ and $$\Vert f_n \Vert_{L^{2}\left( Z^q\right) }^{4} = \Vert f_n \otimes f_n \Vert_{L^{2}\left( Z^{2q}\right) }^{2}= \alpha^4,$$ yields the following expression for the fourth moment of $I^{\hat{N}}_q(f_n)$:
\begin{equation}
\label{expression4mom}
\varphi\left( I^{\hat{N}}_q(f_n)^4\right)  = 2\alpha^4 + \sum_{k=1}^{q-1} \Vert f_n \cont{k} f_n\Vert_{L^{2}\left( Z^{2q-2k}\right) }^{2} +  \sum_{k=1}^{q} \Vert f_n \star_{k}^{k-1} f_n\Vert_{L^{2}\left( Z^{2q-2k+1}\right) }^{2}
\end{equation}
from which the desired conclusion follows. It remains to show that (iv) $\Longrightarrow$ (i). This will be achieved by showing that all the free cumulants of $I^{\hat{N}}_q(f_n)$ of order $\geq 3$ converge to zero as $n$ goes to infinity, hence proving that $I^{\hat{N}}_q(f_n)$ converges to a centered semicircular distribution with variance $\alpha^2$. To do this, we will use the diagram formula \eqref{e:dp} proved in the previous section. Recall that, for every $m\geq 3$,
\begin{equation}
\label{cumulantmdefn}
\kappa_{m}\left( I_{q}^{\hat{N}}(f_n)\right)  = \sum_{\sigma \in {NC}_{\geq 2}\left(\left[mq \right], \pi^{*}  \right) }\int_{Z^{\vert \sigma \vert}}(f_n)_{\sigma}d\mu^{\vert \sigma \vert}, 
\end{equation}
where $(f_n)_{\sigma}$ is defined according to Definition \ref{d:tensor}: the idea is now, for every $\sigma$, to decompose the kernel $(f_n)_{\sigma}$ into two parts, in such a way that the asymptotic behavior of the RHS of \eqref{cumulantmdefn} can be suitably controlled by exploiting (iv) and tameness. To do this, fix $m \geq 3$ and  let $\sigma \in {NC}_{\geq 2}\left(\left[mq \right] , \pi^{*} \right)$. As in Definition \ref{respectfulpartitions}, we denote by $B_1,\ldots,B_m$ the blocks of $\pi^*$ and by $b_j^{(i)}$ the elements of the block $B_i \in \pi^*$, for all $1 \leq j \leq q$ and $1 \leq i \leq m$. By Lemma \ref{l:1} and Remark \ref{remarkdefalternnoncrossing}, we know that $\sigma$ has at least one block $V$ of the form $\left\lbrace b_{q}^{(i)}, b_{1}^{(i+1)} \right\rbrace $ (this specific form for $V$ comes from the fact that $\sigma$ respects $\pi^{*}$, see Remark \ref{remarkdefalternnoncrossing}) for some $1 \leq i \leq q-1$. Note that, using the terminology introduced above, one has that the block $V$ links $B_i$ and $B_{i+1}$. Now denote by $r$ the number of additional blocks of $\sigma$ of size $2$ (other than $V$) linking $B_i$ and $B_{i+1}$. If  $r=0$, then there is only one block of size two (namely $V$) linking $B_i$ and $B_{i+1}$. If $r \geq 1$, then, because $\sigma$ is non--crossing and does not contain any singleton, such $r$ blocks of size two are necessarily given by $\left\lbrace b_{q-1}^{(i)}, b_{2}^{(i+1)} \right\rbrace ,\ldots,\left\lbrace b_{q-r}^{(i)}, b_{r+1}^{(i+1)}\right\rbrace $. Since $\pi^*\vee \sigma = \hat{1}$, one has also that $r+1<q$. Fig. 4 gives a visual representation of what the $r$ blocks of size $2$ look like.
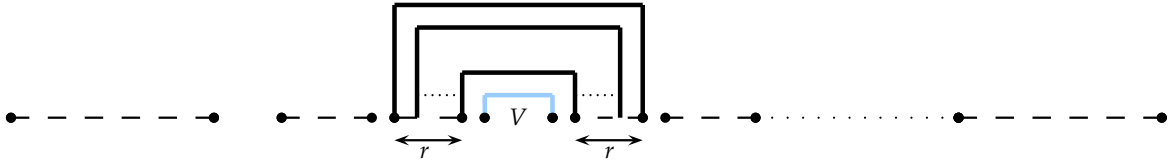
\begin{figure}[H]
\begin{center}
\scalebox{1.5} 
{
\begin{pspicture}(0,-0.73)(10.31,0.71)
\definecolor{color3791}{rgb}{0.00392156862745098,0.00392156862745098,0.00392156862745098}
\definecolor{color5953}{rgb}{0.6,0.8,1}
\definecolor{color6125}{rgb}{0.00392156862745098,0.00392156862745098,0.00392156862745098}
\psdots[dotsize=0.09](0.045,-0.31)
\psdots[dotsize=0.09](1.845,-0.31)
\psdots[dotsize=0.09](2.445,-0.31)
\psdots[dotsize=0.09](8.445,-0.31)
\psdots[dotsize=0.09](4.845,-0.31)
\psdots[dotsize=0.09](6.645,-0.31)
\psdots[dotsize=0.09](8.445,-0.31)
\psdots[dotsize=0.09](10.245,-0.31)
\psdots[dotsize=0.09,linecolor=color3791](4.245,-0.31)
\psdots[dotsize=0.09,linecolor=color3791](4.845,-0.31)
\psdots[dotsize=0.09,linecolor=color3791](6.645,-0.31)
\psdots[dotsize=0.09,linecolor=color3791](8.445,-0.31)
\psdots[dotsize=0.09,linecolor=color3791](1.845,-0.31)
\psdots[dotsize=0.09,linecolor=color3791](2.445,-0.31)
\psdots[dotsize=0.09,linecolor=color3791](0.045,-0.31)
\psdots[dotsize=0.09,linecolor=color3791](10.245,-0.31)
\psdots[dotsize=0.09,linecolor=color3791](4.045,-0.31)
\psdots[dotsize=0.09,linecolor=color3791](8.445,-0.31)
\psdots[dotsize=0.09,linecolor=color3791](3.245,-0.31)
\psdots[dotsize=0.09,linecolor=color3791](3.445,-0.31)
\psline[linewidth=0.02cm,linecolor=color3791,linestyle=dashed,dash=0.16cm 0.16cm](2.445,-0.31)(3.245,-0.31)
\psline[linewidth=0.02cm,linecolor=color3791,linestyle=dashed,dash=0.16cm 0.16cm](3.445,-0.31)(4.045,-0.31)
\usefont{T1}{ppl}{m}{it}
\rput(3.7120314,-0.62){\tiny r}
\psline[linewidth=0.02cm,linecolor=color3791,linestyle=dashed,dash=0.16cm 0.16cm](8.445,-0.31)(10.245,-0.31)
\psdots[dotsize=0.09](5.045,-0.31)
\psdots[dotsize=0.09,linecolor=color3791](5.845,-0.31)
\psdots[dotsize=0.09,linecolor=color3791](5.645,-0.31)
\psline[linewidth=0.02cm,linecolor=color3791,linestyle=dashed,dash=0.16cm 0.16cm](5.245,-0.31)(5.645,-0.31)
\psline[linewidth=0.02cm,linecolor=color3791,linestyle=dashed,dash=0.16cm 0.16cm](5.845,-0.31)(6.645,-0.31)
\psline[linewidth=0.04cm,linecolor=color5953](4.245,-0.11)(4.245,-0.31)
\psline[linewidth=0.04cm,linecolor=color5953](4.245,-0.11)(4.845,-0.11)
\psline[linewidth=0.04cm,linecolor=color5953](4.845,-0.11)(4.845,-0.31)
\psdots[dotsize=0.09,linecolor=color3791](4.245,-0.31)
\psdots[dotsize=0.09,linecolor=color3791](4.845,-0.31)
\psline[linewidth=0.02cm,linecolor=color3791,linestyle=dotted,dotsep=0.12cm](6.645,-0.31)(8.445,-0.31)
\psline[linewidth=0.04cm,linecolor=color6125](4.045,-0.31)(4.045,0.09)
\psline[linewidth=0.04cm,linecolor=color6125](4.045,0.09)(5.045,0.09)
\psline[linewidth=0.04cm,linecolor=color6125](5.045,0.09)(5.045,-0.31)
\psline[linewidth=0.04cm,linecolor=color6125](3.445,-0.31)(3.445,0.69)
\psline[linewidth=0.04cm,linecolor=color6125](3.445,0.69)(5.645,0.69)
\psline[linewidth=0.04cm,linecolor=color6125](5.645,0.69)(5.645,-0.31)
\psline[linewidth=0.04cm,linecolor=color6125](3.645,0.49)(3.645,-0.31)
\psline[linewidth=0.04cm,linecolor=color6125](3.645,0.49)(5.445,0.49)
\psline[linewidth=0.04cm,linecolor=color6125](5.445,0.49)(5.445,-0.31)
\psline[linewidth=0.02cm,linecolor=color3791,linestyle=dotted,dotsep=0.04cm](3.645,-0.11)(4.045,-0.11)
\psline[linewidth=0.02cm,linecolor=color3791,linestyle=dotted,dotsep=0.04cm](5.045,-0.11)(5.445,-0.11)
\psline[linewidth=0.02cm,linecolor=color3791,arrowsize=0.05291667cm 2.0,arrowlength=1.4,arrowinset=0.4]{<->}(3.445,-0.51)(4.045,-0.51)
\psline[linewidth=0.02cm,linecolor=color3791,arrowsize=0.05291667cm 2.0,arrowlength=1.4,arrowinset=0.4]{<->}(5.045,-0.51)(5.645,-0.51)
\usefont{T1}{ppl}{m}{it}
\rput(5.352031,-0.62){\tiny r}
\usefont{T1}{ppl}{m}{it}
\rput(4.538125,-0.3){\tiny V}
\psline[linewidth=0.02cm,linecolor=color3791,linestyle=dashed,dash=0.16cm 0.16cm](0.045,-0.31)(1.845,-0.31)
\psdots[dotsize=0.09,linecolor=color3791](3.445,-0.31)
\psdots[dotsize=0.09,linecolor=color3791](4.045,-0.31)
\psdots[dotsize=0.09,linecolor=color3791](5.045,-0.31)
\psdots[dotsize=0.09,linecolor=color3791](5.645,-0.31)
\end{pspicture} 
}
\caption{\sl Placement of the $r$ additional blocks of size two linking $B_i$ and $B_{i+1}$.}
\end{center}
\end{figure}
\noindent There are two additional categories of blocks of $\sigma$ that are allowed to link the remaining unassigned elements of $B_i$ and $B_{i+1}$: (a) blocks composed of an element of $B_i$, an element of $B_{i+1}$, and other elements from $[mq]\backslash (B_{i}\cup B_{i+1})$; (b) blocks composed of one element of $B_{i}$ (resp. $B_{i+1}$), no elements of $B_{i+1}$ (resp. $B_{i}$), and elements from $[mq]\backslash (B_{i}\cup B_{i+1})$. We denote the number of blocks of type (a) by $p$. It is immediately seen that $p=0$ or $1$, in view of the non--crossing nature of $\sigma$. Also, the number of blocks of $\sigma$ linking $B_{i}$ to other blocks and not to $B_{i+1}$ is the same as the amount of blocks of $\sigma$ linking $B_{i+1}$ to other blocks and not to $B_{i}$.


\noindent Now, for integers $r,p$ as above, define the following function $f_{i,n}$, in $2q-r-1-p$ variables:
\begin{eqnarray*}
f_{i,n}\left(t_1,\ldots,t_{q-r-1-p},\gamma_{p},z_{1},\ldots, z_{r+1},s_1,\ldots,s_{q-r-1-p} \right) &=& f_n\left(t_{q-r-1-p},\ldots,t_1,\gamma_{p},z_{r+1},\ldots,z_1\right)\\
&& f_n\left(z_1,\ldots,z_{r+1},\gamma_{p},s_{1},\ldots,s_{q-r-1-p} \right),
\end{eqnarray*}
with the convention that the variable $\gamma_{p}$ is deleted if $p=0$. Observe that 
\begin{eqnarray*}
&& \int_{Z^{r+1}}f_{i,n}\left(t_1,\ldots,t_{q-r-1-p},\gamma_{p},z_{1},\ldots, z_{r+1},s_1,\ldots,s_{q-r-1-p} \right)\mu\left( dz_{1}\right) \cdots \mu\left( dz_{r+1}\right)  \\
&& \qquad\qquad\qquad\qquad\qquad\qquad\qquad = f_n \star_{r+1+p}^{r+1} f_n\left(t_1,\ldots,t_{q-r-1-p},\gamma_{p},s_1,\ldots,s_{q-r-1-p} \right),
\end{eqnarray*}
and consequently
\begin{eqnarray*}
\label{decompfsig}
\int_{Z^{\vert \sigma \vert}}(f_n)_{\sigma}d\mu^{\vert \sigma \vert} &=& \int_{Z^{\vert \sigma \vert - r -1}}f_n \star_{r+1+p}^{r+1} f_n \times \mathscr{G}_n^{\sigma, r} d\mu^{\vert \sigma\vert - r -1}\\
&=& \int_{Z^{2q - 2r -2 -p}}f_n \star_{r+1+p}^{r+1} f_n \left( \int_{Z^{\vert \sigma\vert +r+1+p-2q}}\mathscr{G}_n^{\sigma, r} d\mu^{\vert \sigma\vert +r+1+p-2q}\right) d\mu^{2q - 2r -2 -p},
\end{eqnarray*}
where $\mathscr{G}_n^{\sigma, r}$ is a function whose argument is composed of the $\vert\sigma\vert - r -1$ variables that are not integrated out in the definition of the star contraction, and the integral between brackets  in the second equality is realized by integrating with respect to all variables that are not in the argument of $f_n \star_{r+1+p}^{r+1} f_n$ (if there are no such variables, then the second equality is immaterial). Replacing in \eqref{cumulantmdefn}  and applying the Cauchy--Schwarz inequality yields
\begin{eqnarray*}
 \left| \kappa_{m}\left( I_{q}^{\hat{N}}(f_n)\right) \right| &\leq & \sum_{\sigma \in {NC}_{\geq 2}\left(\left[mq \right], \pi^{*}  \right) }  \left\| f_n \star_{r+1+p}^{r+1} f_n \right\|_{L^2\left(Z^{2q - 2r -2 -p} \right) } \\
&& \times  \sqrt{\int_{Z^{2q - 2r -2 -p}} \left|\int_{Z^{\vert \sigma\vert +r+1+p-2q}}\mathscr{G}_n^{\sigma, r} d\mu^{\vert \sigma\vert +r+1+p-2q}\right|^2 d\mu^{2q - 2r -2 -p}}.
\end{eqnarray*}
Now observe that, since $\{f_n\}$ is tamed, the sequence $$\left\lbrace \int_{Z^{2q - 2r -2 -p}} \left| \int_{Z^{\vert \sigma\vert +r+1+p-2q}}\mathscr{G}_n^{\sigma, r} d\mu^{\vert \sigma\vert +r+1+p-2q}\right|^2 d\mu^{2q - 2r -2 -p} \colon n\geq 1\right\rbrace $$ is bounded. Recalling that $p \in \left\lbrace 0,1 \right\rbrace $ and using the assumption that, as $n \rightarrow \infty$, $\| f_n \cont{k} f_n\|_{L^2(Z^{2q-2k})} \rightarrow 0$ for all $k \in \left\lbrace 1,\ldots,q-1 \right\rbrace $ and $\left\| f_n \star_{k}^{k-1} f_n \right\|_{L^2(Z^{2q-2k+1})} \rightarrow 0$ for all $k \in \left\lbrace 2,\ldots,q \right\rbrace$, we finally obtain that, as $n \rightarrow \infty$, $$\kappa_{m}\left( I_{q}^{\hat{N}}(f_n)\right) \longrightarrow 0,$$ thus concluding the proof.\qed
\\~\\ The following consequence of Theorem \ref{fourthcumulant}, which is an important generalization of Proposition \ref{p:noc1}--(i), establishes the fact that no non--trivial self--adjoint multiple integral in $ L_b(\mathcal{S}(\hat{N}),\varphi)$ is distributed according to the semicircular law.
\begin{corollary}
\label{corollarynosemiinpoissonchaos}
Fix $q\geq 2$, and consider a non--zero bounded mirror--symmetric function $f$ with bounded support in $L^2(Z^q)$. Then, the free Poisson multiple integral $I_{q}^{\hat{N}}\left( f\right) $ satisfies $\varphi\left(I_{q}^{\hat{N}}\left( f\right)^4 \right) > 2\varphi\left(I_{q}^{\hat{N}}\left( f\right)^2 \right)^2 $. In particular, the distribution of $I_{q}^{\hat{N}}\left( f\right)$ cannot be semicircular.
\end{corollary}
\noindent{\bf Proof}. 
By rescaling, we may assume that $\Vert f \Vert_{L^2\left( Z^{q}\right)}^{2} = 1$; in this case,
equation \eqref{expression4mom} shows that $\varphi\left(I_{q}^{\hat{N}}\left( f\right)^4 \right) \geq 2\varphi\left(I_{q}^{\hat{N}}\left( f\right)^2 \right)^2 $. To achieve a contradiction, we assume that $\varphi\left(I_{q}^{\hat{N}}\left( f\right)^4 \right) = 2\varphi\left(I_{q}^{\hat{N}}\left( f\right)^2 \right)^2 = 2$, which would be the case if $I_{q}^{\hat{N}}\left( f\right)$ had a semicircular law. According to Theorem \ref{fourthcumulant}, this yields that  $\Vert f \cont{q-1} f \Vert_{L^2\left( Z^{2}\right)}^2 = 0$. By applying the same argument as in \cite[Proof of Corollary 1.7]{knps}, we eventually infer that $f=0$, $\mu^q$--a.e., which contradicts the fact that $\Vert f \Vert_{L^2\left( Z^{q}\right)}^{2} = 1$. \qed

\section{Applications to transfer principles}\label{s:cont}

\noindent As before, we work on the measure space $(Z,\mathscr{Z},\mu) = \left(\R^d, \mathscr{B}(\R^d), \mu\right)$, where $\mu$ stands for the Lebesgue measure. We write $\mathscr{Z}_{\mu} = \{ B\in \mathscr{Z}: \mu(B)< \infty \}$. The notation
$\eta = \{\eta(B) : B\in \mathscr{Z}_{\mu} \} $ is used to indicate a classical {\it Poisson measure} on $\left(Z, \mathscr{Z}\right)$ with control measure $\mu$. This means that $\eta $ is a collection of random variables defined on some probability space $(\Omega, \mathscr{F}, \mathbb{P}) $, indexed by
the elements of $\mathscr{Z}_{\mu} $ and such that: (i) for every $B,C \in \mathscr{Z}_{\mu}$ such that $B \cap C = \varnothing$, the random variables $ \eta(B)$ and $ \eta(C)$ are independent;  (ii) for every $B \in \mathscr{Z}_{\mu} $, $\eta(B)$ has a Poisson distribution with mean $\mu(B)$. We shall also write $$\hat{\eta}(B) = \eta(B) - \mu(B), \quad B\in \mathscr{Z}_{\mu},$$ and $\hat{\eta} = \{\hat{\eta}(B) : B\in \mathscr{Z}_{\mu} \}$. Similarly, we denote by $W$ a centered random Gaussian measure (with control measure $\mu$) on $\left(Z, \mathscr{Z}\right)$. We will write $I^{\hat{\eta}}$ and $I^{W}$ to denote multiple stochastic integrals with respect to $\hat{\eta}$ and $W$, respectively. For more information about the construction of classical stochastic integrals with respect to $\hat{\eta}$ and $W$, see the monograph \cite{PecTaq}.

\medskip

\noindent We start by recording an easy consequence of \cite[Theorem 3.1]{LacPec}, showing that a fourth moment theorem holds also in the classical case whenever one integrates kernels with a constant sign, and also that a partial (i.e. `one-directional') transfer principle takes place, connecting normal and semicircular limits, respectively in the classical and free setting. 

\begin{prop}\label{p:z} For $q\geq 1$, let $\{f_n : n\geq 1\}\subset L^2(Z^q)$ be sequence of real--valued, tamed and fully symmetric kernels, verifying the normalization condition $q!\|f_n\|^2_{L^2(Z^q)} \to 1$.
\begin{itemize}

\item[\rm (A)] If $q=1$, then the following three conditions are equivalent, as $n\to\infty$: {\rm (1)} $I^{\hat{\eta}}_1(f_n)$ converges in law to $\mathcal{N}(0,1)$, {\rm  (2)} $\mathds{E}[I^{\hat{\eta}}_1(f_n)^4]\to 3$, {\rm (3)} $I^{\hat{N}}_1(f_n)$ converges in law to $\mathcal{S}(0,1)$.

\item[\rm (B)] If $q\geq 2$ and $f_n\geq 0$, then the following three conditions are equivalent, as $n\to\infty$: {\rm (1)} $I^{\hat{\eta}}_q(f_n)$ converges in law to $\mathcal{N}(0,1)$, {\rm (2)} $\mathds{E}[I^{\hat{\eta}}_q(f_n)^4]\to 3$, {\rm (3)} $\| f_n\|_{L^4(Z^q)}\to 0$ and $\|f_n\star_r^l f_n\|_{L^2(Z^{2q-r-l})}\to 0$, for every $r=1,...,q$ and every $l=1,...,r\wedge (q-1)$. Moreover, if any of conditions {\rm (1)--(3)} holds, one has that $I^{\hat{N}}_q(f_n)$ converges in law to $\mathcal{S}(0,q!^{-1})$.
\end{itemize}
\end{prop}
\noindent{\bf Proof.} The equivalence between (1) and (2) in Part (A) follows from tameness and tom the fact that, for $m\geq 3$ the classical cumulant of $I^{\hat{\eta}}_1(f_n)$ equals $\int_Z f_n^m d\mu$ (see \cite[Section 7]{PecTaq}). The equivalence with Point (3) is an easy consequence of Proposition \ref{p:noc1}. The first part of Point (B) is a consequence of \cite[Theorem 3.1]{LacPec}, whereas the final assertion follows from Theorem \ref{fourthcumulant}, and from the fact that, by definition, $f_n\star_r^r f_n =  f_n \cont{r} f_n$.
\qed

\medskip

\noindent As announced in Theorem \ref{t:bp}--(B) of the Introduction, we will devote the last two sections of the paper to the construction of an explicit collection of counterexamples, showing that a full version of the transfer principle in Part (B) fails for every order of $q>1$ of integration, that is: it is {\it not} true that, if $I^{\hat{N}}_q(f_n)$ converges to a semicircular limit, then $I^{\hat{\eta}}_q(f_n)$ must verify a CLT. The key point of our construction consists in producing, for every $q\geq 2$, a sequence of tamed positive symmetric kernels $\{f_n\}\subset L^2(Z^q)$ such that $\|f_n\star_r^l f_n\|_{L^2(Z^{2q-r-l})}\to 0$, $\forall r=1,...,q$ and $\forall l=1,...,r\wedge (q-1)$, and such that $\| f_n\|_{L^4(Z^q)}$ converges to some strictly positive limit, in such a way that one cannot apply the previous Proposition \ref{p:z}--(A). Verifying that the actual limit of the sequence $I^{\hat{\eta}}_q(f_n)$ is a centered Poisson distribution will be a relatively simple combinatorial task.

\subsection{A sequence of kernels}
\noindent Fix an integer $q \geq 2$, and let $\{r_n : n\geq 1\}$ be a sequence of positive real numbers such that
\begin{equation}
\label{defrn}
r_n^d \underset{n \rightarrow \infty}{\sim} n^{-\frac{1}{q-1}}
\end{equation}
For every $n\geq 1$, we write $\Q_n = \left[-\frac{1}{2}n^{1/d},\frac{1}{2}n^{1/d} \right]^{d}$. The collection of symmetric kernels $\left\lbrace f_n \colon n \geq 1 \right\rbrace \subset L^2\left(Z^q \right)$ that are of interest in the present section are defined by 
\begin{equation}
\label{noyauxfn}
f_n(x_1,\ldots,x_q) := \frac{1}{q!}\mathds{1}_{\Q^q_n}(x_1,...,x_q)\times \mathds{1}_{\{ 0< \left\| x_i -x_j\right\| \leq r_n , \forall\  1\leq i,j\leq q \} },
\end{equation}
that is, the function $f_n(x_1,\ldots,x_q)$ is defined to be $1/q!$ if all of its arguments lie inside $\Q_n$ and are at a distance of at least $r_n$, and vanishes otherwise. Applying Lemma \ref{l:tamed} in the case $z_n = 2^{-1} n^{1/d}$, $M_n = 1/q!$ and $\alpha_n = r_n$ yields immediately that such a sequence of kernels is tamed.

\begin{remark}{\rm The kernels $f_n$ admit a straightforward geometric interpretation. Let $V$ be a finite subset of $\Q_n$, and denote by $G_n = (V,E)$ the undirected graph obtained by connecting two distinct points $v_1,v_2 \in V$ if and only if their Euclidean distance is less than $r_n$. Then, the quantity
\[
K_n := \sum_{(v_1,...v_q)\in V^q} f_n(v_1,...,v_q)
\]
equals the number of subsets $\{ v_1,...,v_q\}\subset V$ (of size exactly $q$) forming a {\it clique}, that is, such that the restriction of $G_n$ to $\{v_1,...,v_q\}$ is a complete graph. Note that one could build many more sequences of kernels having the same asymptotic properties as those of the kernels $f_n$ (and therefore violating the transfer principle), by simply replacing the complete graph with any connected graph with $q$ vertices.}
\end{remark}

\noindent The following lemma provides insights into the behavior of the sequence $\left\lbrace f_n \colon n \geq 1 \right\rbrace$ and will play a central role in the sequel. The simple proof involves computations that are completely analogous to the ones leading to the proof of Lemma \ref{l:tamed}, and are therefore left to the reader.
\begin{lemma}
\label{lemmacontnoayuxfn}
Let the sequence $\left\lbrace f_n \colon n \geq 1 \right\rbrace$ be defined as in \eqref{noyauxfn}. Then, the following facts hold.
\begin{enumerate}
\item[\rm (i)] There exist a constant $\alpha > 0$ such that $\alpha = \underset{n \rightarrow \infty}{\rm lim} q!\Vert f_n \Vert_{L^2\left( Z^q\right) }^{2}$.
\item[\rm (ii)] For all $r = 1,\ldots ,q$ and $l = 1,\ldots , r\wedge (q-1)$, $\Vert f_n \star_{r}^{l} f_n \Vert_{L^2\left( Z^{2q - r - l}\right) }^{2} \underset{n \rightarrow \infty}{\longrightarrow} 0$.
\end{enumerate}
\end{lemma}

\subsection{Counterexamples to the transfer principle}
The next statement characterizes the asymptotic behavior of the multiple integrals of the sequence $\{f_n\}$, realized with respect to several different (classical and free) measures. Given $\alpha>0$, we will denote by $\hat{\rm{P}}\rm{o}\left(\alpha\right)$ a classical centered Poisson distribution of parameter $\alpha$ (that is, $Y\sim \hat{\rm{P}}\rm{o}\left(\alpha\right)$ if and only if $Y +\alpha$ has a usual Poisson distribution with parameter $\alpha$). Recall that the distribution of $\hat{\rm{P}}\rm{o}\left(\alpha\right)$ is characterized by the fact of having all classical cumulants of order $\geq 2$ equal to $\alpha$ (see \cite[Chapter 3]{PecTaq}).
\begin{prop}
\label{remarksontransferanduni}
Fix an integer $q\geq 2$, and consider the sequence of functions $\left\lbrace f_n \colon n \geq 1 \right\rbrace $ defined in \eqref{noyauxfn}. Let $\alpha$ be the constant defined in Lemma \ref{lemmacontnoayuxfn}. Then, as $n \rightarrow \infty$, the following convergences in law are in order:
\begin{enumerate}
\item[\rm (i)] $I_{q}^{W}\left( f_n\right) \longrightarrow \mathcal{N}(0,\alpha)$;
\item[\rm (ii)] $I_{q}^{\hat{\eta}}\left( f_n\right) \longrightarrow \hat{\rm{P}}\rm{o}\left(\alpha\right)$;
\item[\rm (iii)] $I_{q}^{S}\left( f_n\right) \longrightarrow \mathcal{S}(0,\alpha)$;
\item[\rm (iv)] $I_{q}^{\hat{N}}\left( f_n\right) \longrightarrow \mathcal{S}(0,\alpha)$.
\end{enumerate}
\end{prop}
\noindent{\bf Proof}. In order to prove Point (i), it suffices to combine Lemma \ref{lemmacontnoayuxfn} with the fourth moment theorem proved in \cite[Theorem 1]{nuapec1}. Applying the transfer principle in \cite[Theorem 1.8]{knps} yields Point (iii), whereas Point (iv) follows directly from the combination of Lemma \ref{lemmacontnoayuxfn} and Theorem \ref{fourthcumulant}. It remains to prove point (ii). To do this, we shall apply the method of moments, and prove that, for every integer $m\geq 3$,
\[
\chi_m(I_{q}^{\hat{\eta}}\left( f_n\right))\to \alpha,
\]
where the symbol $\chi_m(Y)$ denotes the $m^{\mbox{\tiny{th}}}$ classical cumulant of a given random variable $Y$, as defined e.g. in \cite[Chapter 3]{PecTaq}. Now fix $m\geq 3$ and write $\pi^*$ in order to indicate (as before) the partition of $[mq]$ composed of $m$ consecutive blocks of size $q$. According to the refinement of the content of \cite[Section 7.4]{PecTaq} recently proved in \cite[Theorem 3.4]{LPST}, one has that the following diagram formula holds
\[
\chi_{m}\left( I_{q}^{\hat{\eta}}(f)\right)  = \sum_{\sigma}\int_{Z^{\vert \sigma \vert}}f_{\sigma}d\mu^{\vert \sigma \vert},
\]
where the sum runs over all partitions $\sigma\in \mathcal{P}([mq])$ such that $\sigma\vee \pi^* =\hat{1}$ and $\sigma\wedge \pi^* =\hat{0}$. Reasoning as in the proof of Theorem \ref{fourthcumulant}, one sees that there are two possibilities: either (a) $|\sigma| >q$, and 
\[
\left| \int_{Z^{\vert \sigma \vert}}f_{\sigma}d\mu^{\vert \sigma \vert}\right| = O(1)\,  \Vert f_n \star_{r}^{l} f_n \Vert_{L^2\left( Z^{2q - r - l}\right) }\to 0,
\]
where $O(1)$ stands for a bounded numerical sequence, $r\in\{1,...,q\}$ and $l\in \{1,...,r\wedge (q-1)\}$, or (b) $\sigma$ has exactly $q$ blocks of size $m$, and 
\[
 \int_{Z^{\vert \sigma \vert}}f_{\sigma}d\mu^{\vert \sigma \vert} \to \frac{1}{q!^{m-1}} \alpha.
\]
Since there are exactly $q!^{m-1}$ partitions as in (b), we deduce immediately the desired conclusion.
\qed
\\~\\

\begin{remark} 
{\rm \textcolor{black}{\begin{enumerate}
\item[(a)] The fact that Points (i) and (iii) in the previous statement hold simultaneously (they are in fact equivalent) is a demonstration of the transfer principle between the Wiener and Wigner chaos proved in \cite[Theorem 1.8]{knps}, and recalled in Theorem \ref{t:knps}--(B).
\item[(b)] In \cite{npr}, the authors proved a new universality phenomenon for \textit{homogeneous sums} of classical random variables. Roughly speaking, for $q\geq 2$ this result implies that, if $\{f_n : n\geq 1\}$ is a normalized sequence of kernels in the class $\mathscr{E}_q^{00}$ (that is, the collection of all simple kernels vanishing on diagonals), then the convergence in law $I_{q}^{W}\left( f_n\right) \longrightarrow \mathcal{N}(0,1)$ takes place if and only if $I_{q}^{\hat{\eta}}\left( f_n\right) \longrightarrow \mathcal{N}(0,1)$. It is natural to ask whether this phenomenon extends to more general sequences of non--simple kernels: the fact that Point (i) and (ii) of Proposition \ref{remarksontransferanduni} hold shows that the answer is negative for every order $q$ of integration.
\end{enumerate}
}}
\end{remark}

\section{Two proofs}\label{s:proofs}

\subsection{Proof of Theorem \ref{t:mm}} By definition, the multiple integrals appearing on both sides of the equalities \eqref{e:mulp}--\eqref{e:muls} are elements of the unit algebras generated, respectively, by $\hat{N}$ and $S$. It is therefore enough to prove these relations in the special case where the free Poisson and semicircular measures are both defined on the free Fock space associated with $L^2(Z)$. For the convenience of the reader, we include a quick discussion of these concepts.

\medskip

\noindent{\it Preliminary definitions}. Let $\HH = L^2(Z)$, and let $\mathscr{F}_{0}\left(\HH \right)$ be its algebraic Fock space, defined by \[\mathscr{F}_{0}\left(\HH \right) = \bigoplus_{n=0}^{\infty}\HH^{\otimes n},\] where the direct sum and tensor products are Hilbert space operations, and $\HH^{\otimes 0}\equiv \mathbb{C}\Omega$ is defined to be a one dimensional complex space with a distinguished unit basis vector $\Omega$, called the {\it vacuum vector}. For $h \in \HH$, we define the {\it creation} operator $a^{+}(h)$, the {\it annihilation} operator $a^{-}(h)$ and the {\it gauge} operator $a^{o}(h)$ on $\mathscr{F}_{0}\left(\HH \right)$  in terms of their actions on $\Omega$ and on $n$--tensors:
\begin{eqnarray*}
&& a^{+}(h)\Omega = h, \quad a^{+}(h)(f_1 \otimes f_2 \otimes \cdots \otimes f_n) = h\otimes f_1 \otimes f_2 \otimes \cdots \otimes f_n,  \\
&& a^{-}(h)\Omega = 0, \quad a^{-}(h)(f_1 \otimes f_2 \otimes \cdots \otimes f_n) = \left\langle h,f_1 \right\rangle _{\HH} f_2 \otimes \cdots \otimes f_n,  \\
&& a^{o}(h)\Omega = 0, \quad a^{o}(h)(f_1 \otimes f_2 \otimes \cdots \otimes f_n) = (hf_1) \otimes f_2 \otimes \cdots \otimes f_n.
\end{eqnarray*}

\medskip

\noindent{\it Poisson case.} Let $\hat{N}_0$ be free random measure realized on the Fock space as $a^+ + a^- + a^o$, i.e. for every $A\in \mathcal{Z}_\mu$, $\hat{N}_0(A)$ is realized as $a^+\left( \mathds{1}_{A}\right)  + a^-\left( \mathds{1}_{A}\right) + a^o\left( \mathds{1}_{A}\right)$. Denote by $\mathcal{S}(\hat{N}_0)$ the unital algebra generated by all operators $\hat{N}_0(A)$. Then, the following facts are well--known (see e.g. \cite{ans0, ans1, speicher90}): 
\begin{itemize}

\item[--] The mapping $\varphi : \mathcal{S}(\hat{N}_0) \to \mathbb{C} : T\mapsto \langle T\Omega, \Omega\rangle_{\mathscr{F}_{0}\left(\HH \right)}$ is a faithful normal trace on $\mathcal{S}(\hat{N}_0)$.

\item[--] The operators $\{\hat{N}_0(A) : A\in \mathcal{Z}_\mu\}$ form a free Poisson measure with control $\mu$.

\item[--] The mapping $\mathcal{S}(\hat{N}_0) \to\mathscr{F}_{0}\left(\HH \right) :  T\mapsto T\Omega$ is an injective isometry.

\item[--] Any multiple integral of the form $I^{\hat{N}}_q(f)$, as defined in formula (\ref{e:1}), is characterized by the relation $I^{\hat{N}_0}_q(f)\Omega = f_1\otimes \cdots \otimes f_q$.
\end{itemize}
In view of Lemma \ref{l:taborn}, of the linearity of multiple integrals and of the isometric property stated in Proposition \ref{p:iso}, it is enough to prove formula \eqref{e:mulp} in the case where $g = g_1\otimes \cdots \otimes g_n\in \mathscr{E}_n$, and $f\in \mathscr{E}_m^0$ has the form $f = f_1\otimes \cdots \otimes f_m$, where $f_j = \mathds{1}_{A_j}$ and $A_j \cap A_{j+1} = \emptyset$. One has that
\begin{eqnarray*}
&&I_m^{\hat{N}_0}(f_1\otimes \cdots \otimes f_m)I_n^{\hat{N}_0}(g_1\otimes \cdots \otimes g_n)\Omega =I_m^{\hat{N}_0}(f_1\otimes \cdots \otimes f_m)\, g_1\otimes \cdots \otimes g_n\\
&&= [a^+(f_1)+a^-(f_1)+a^o(f_1)]\cdots [a^+(f_m)+a^-(f_m)+a^o(f_m)] \, g_1\otimes \cdots \otimes g_n.
\end{eqnarray*}
The last expression can be written as a sum of $3^m$ terms of the type 
\begin{equation}\label{e:ggg}
e_m(f_1) e_{m-1}(f_2) \cdots e_1(f_m)  ( g_1\otimes \cdots \otimes g_n),
\end{equation}
were $e_j$ equals either one of the symbols $a^+,\, a^-$ or $a^o$, for $j=1,\ldots,m$. We identify each of these summands with the string $e_me_{m-1}\cdots e_1$: note that such a string is labeled in increasing order from the right to the left, in agreement with the fact that the non--commutative setting requires one to keep track of the order in which the operators act. Since $A_{j-1}\cap A_{j} = \emptyset$ for every $j=2,\ldots,m$, it is clear that the only non--vanishing terms in the sum are those corresponding to strings $e_me_{m-1}\cdots e_1$ that obey the following rules for every $j=1,\ldots,m-1$: (i) if $e_j= a^+$ or $e_j= a^o$, then $e_{j+1} = a^+$, (ii)  if $e_j = a^-$, then $e_{j+1}$ can be either $a^-,\, a^+$ or $a^o$, except in the case where $j+1>n$ (which only can hold whenever $m>n$), where one must have $e_{j+1} = a^+$. We can explicitly describe all strings satisfying properties (i)--(ii): (a) the trivial string such that $e_i = a^+$ for every $i=1,\ldots,m$, for which the expression in (\ref{e:ggg}) equals $f\cont{0} g$, (b) strings such that $e_i = a^-$ for all $i=1,\ldots ,k$, for some $k\leq m\wedge n$, and $e_i = a^+$ for $i>k$, for which \eqref{e:ggg} is equal to $f\cont{k} g$, (c) the string for which $e_1 = a^o$ and $e_i=a^+$ for all $i>k$, for which \eqref{e:ggg} is equal to $f\star_1^0 g$, and finally (d) strings for which there exists $k\in \{2,\ldots,m\wedge n\}$ such that $e_i = a^-$ for every $i=1,\ldots,k-1$, $e_k = a^o$ and $e_i=a^+$ for all $i>k$, for which \eqref{e:ggg} equals $f\star_k^{k-1} g$. Summing up, we just proved that
\begin{eqnarray*}
&&I_m^{\hat{N}_0}(f_1\otimes \cdots \otimes f_m)I_n^{\hat{N}_0}(g_1\otimes \cdots \otimes g_n)\Omega = \sum_{k=0}^{m \wedge n} f \cont{k} g  +  \sum_{k=1}^{m \wedge n} f \star_{k}^{k-1} g,
\end{eqnarray*}
which is the desired conclusion.

\medskip

\noindent{\it Semicircular case}. The proof follows exactly the same structure, except that one has now to work with the unital algebra $\mathcal{S}(S_0)$ generated by the operators $S_0(A)$ realized on the Fock space as $ a^-(\mathds{1}_A) + a^+(\mathds{1}_A)$, $A\in \mathcal{Z}_\mu$, endowed with the faithful tracial operator $\varphi(T) = \langle T\Omega, \Omega\rangle_{\mathscr{F}_{0}\left(\HH \right)}$ introduced above. It is well--known (see e.g. \cite{bianespeicher}) that $S_0$ is a semicircular measure with control $\mu$, and also that multiple integrals of the form $I^{S_0}_q(f)$, as given in formula (\ref{e:2}), are characterized by the relation $I^{S_0}_q(f)\Omega = f_1\otimes \cdots \otimes f_q$. Now take $g = g_1\otimes \cdots \otimes g_n\in \mathscr{E}_n$, and $f\in \mathscr{E}_m^0$ with the form $f = f_1\otimes \cdots \otimes f_m$, where $f_j = \mathds{1}_{A_j}$ and $A_j \cap A_{j+1} = \emptyset$. Reasoning as in the first part of the proof, one sees that the quantity $I_m^{S_0}(f_1\otimes \cdots \otimes f_m)I_n^{S_0}(g_1\otimes \cdots \otimes g_n)\Omega$ can be expanded into a sum of $2^m$ terms of the type $e_m(f_1) e_{m-1}(f_2) \cdots e_1(f_m)  ( g_1\otimes \cdots \otimes g_n)$, with $e_j$ equal to either $a^+$ or $a^-$, for $j=1,\ldots,m$. The conclusion follows by applying the same combinatorial arguments exploited above.

\qed

\subsection{Proof of Theorem \ref{freediagramformulae}}
\noindent {\it Preparation.} We start by considering the case $f\in \mathscr{E}^{00}_q$, that is, the kernel $f : Z^q \to \mathbb{C}$ has the form 
\begin{equation}
\label{fonctionssimples}
f\left( t_1,\ldots,t_q\right)  = \sum_{i_1,\ldots,i_q}c_{i_1\cdots i_q}\mathds{1}_{A_{i_1} \times \cdots \times A_{i_q}}\left( t_1,\ldots,t_q\right),
\end{equation}
where the sum is finite, $A_{i_1},\ldots,A_{i_q}\in \mathcal{Z}_{\mu}$ are disjoint and bounded, and the complex coefficients $c_{i_1\cdots i_q}$ are zero if $i_j = i_k$ for some $j\neq k$. By \eqref{e:wickmaps}, for every $f$ of the form (\ref{fonctionssimples}) one has that
$$I_q^{\mathfrak{M}}(f)=\sum_{i_1,\ldots,i_q}c_{i_1\cdots i_q}{\mathfrak{M}}\left( A_{i_1}\right) \cdots {\mathfrak{M}}\left( A_{i_q}\right).$$
We will adopt the notations $A_{q(l-1) +s} =  A_{s}^{(l)}$ and $\Delta_{q(l-1) +s} =  {\mathfrak{M}}\left( A_{q(l-1) +s}\right)$, $l = 1,\ldots, m$, $s = 1,\ldots, q$, to refer to the collection of $mq$ random variables (that we shall sometimes refer to as the `increments' of ${\mathfrak{M}}$ for simplicity) $$\left\lbrace {\mathfrak{M}}\left( A_{1}^{(1)}\right),\ldots, {\mathfrak{M}}\left( A_{q}^{(m)}\right) \right\rbrace .$$ Recall that free cumulants are multilinear functionals, hence yielding the following expression for the free cumulant of order $m$ of $I_{q}^{{\mathfrak{M}}}(f)$:
\begin{equation}
\label{prespeicher}
\kappa_{m}\left( I_{q}^{{\mathfrak{M}}}(f)\right)  = \sum_{i_{1}^{(1)},\ldots,i_{q}^{(1)}}\cdots \sum_{i_{1}^{(m)},\ldots,i_{q}^{(m)}} c_{i_{1}^{(1)}\cdots i_{q}^{(1)}}\cdots c_{i_{1}^{(m)}\cdots i_{q}^{(m)}} \kappa \left(  \prod_{s=1}^{q} \Delta_{s}         ,\ldots, \prod_{s=1}^{q} \Delta_{q(m-1)+s} \right),
\end{equation}
where the joint cumulant notation is taken from \cite[p. 175]{NicSpe}. Using \cite[Theorem 11.12, point 2]{NicSpe}, we can express the cumulants appearing in (\ref{prespeicher}) in the following way:
\begin{equation}
\label{sumsurpartitions}
\kappa \left(  \prod_{s=1}^{q} \Delta_{s},\ldots, \prod_{s=1}^{q} \Delta_{q(m-1)+s} \right) = \sum_{\substack{\sigma \in {NC}(qm )\\\sigma \vee \pi^{*} = \hat{1}}}\kappa_{\sigma}\left[\Delta_{1},\ldots,\Delta_{q(m-1) + q} \right],
\end{equation}
where the quantities $\kappa_{\sigma}\left[\Delta_{1},\ldots,\Delta_{q(m-1) + q} \right]$ are defined by the relation
\begin{equation*}
\kappa_{\sigma}\left[\Delta_{1},\ldots,\Delta_{q(m-1) + q} \right] = \prod_{V \in \sigma}\kappa(V)\left[\Delta_{1},\ldots,\Delta_{q(m-1) + q} \right] 
\end{equation*}
with $\kappa(V)\left[\Delta_{1},\ldots,\Delta_{q(m-1) + q} \right] = \kappa\left(\Delta_{q(l-1) +s},\ \ q(l-1) +s \in V \right)$. In the sum appearing on the right hand side of (\ref{sumsurpartitions}), consider a partition $\sigma$ such that $\sigma \wedge \pi^{*} \neq \hat{0}$. This property implies that $\sigma$ does not respect $\pi^{*}$ and hence that there exists a block $V^{*} \in \sigma$ such that the argument of $\kappa(V^{*})\left[\Delta_{1},\ldots,\Delta_{q(m-1) + q} \right]$ contains a non--empty collection of increments of the measure ${\mathfrak{M}}$ over disjoint sets. In view of the freeness of the increments of the measure ${\mathfrak{M}}$, we infer that $\kappa(V^{*})\left[\Delta_{1},\ldots,\Delta_{q(m-1) + q} \right] = 0$ and hence $\kappa_{\sigma}\left[\Delta_{1},\ldots,\Delta_{q(m-1) + q} \right] = 0$. Also, if there exists a block $V^{'} \in \sigma$ such that $\vert V^{'} \vert = 1$ (assume $V^{'} = \left\lbrace 1 \right\rbrace$ without any loss of generality), then $\kappa(V^{'})\left[\Delta_{1},\ldots,\Delta_{q(m-1) + q} \right] = \kappa \left( \Delta_1 \right)  = \varphi\left( \Delta_1\right) =0$ because ${\mathfrak{M}}$ is centered. Therefore, (\ref{sumsurpartitions}) can be rewritten as 
\begin{equation}
\label{sumsurnc}
\kappa \left(  \prod_{s=1}^{q} \Delta_{s},\ldots, \prod_{s=1}^{q} \Delta_{q(m-1)+s} \right) = \sum_{\sigma \in {NC}_{\geq 2}\left(\left[mq \right], \pi^{*}  \right) }\kappa_{\sigma}\left[\Delta_{1},\ldots,\Delta_{q(m-1) + q} \right].
\end{equation}
Using once again the freeness of the increments of the measure ${\mathfrak{M}}$, it is clear that if the increments appearing in the expression $\kappa(V)\left[\Delta_{1},\ldots,\Delta_{q(m-1) + q} \right]$ are different, this quantity is zero. Consequently, we can write $\kappa(V)\left[\Delta_{1},\ldots,\Delta_{q(m-1) + q} \right] = \kappa_{\vert V \vert}\left(\Delta_{b_1^{V,\sigma}}\right)$, where $b_1^{V,\sigma}$ is the first element of the block $V \in \sigma$. We will now differentiate between the case where ${\mathfrak{M}} = \hat{N}$ and the case where ${\mathfrak{M}} =S$.

\medskip

\noindent {\it Case 1: ${\mathfrak{M}} = \hat{N}$.} As above, we consider a non--diagonal $f\in \mathscr{E}^{00}_q$. We will use the fact that, for any $n \geq 2$ and any $A \in \mathcal{Z}_{\mu}$, $\kappa_{n}\left(\hat{N}\left(A \right) \right) = \mu\left( A\right) $. Consequently, one can rewrite (\ref{sumsurnc}) as
\begin{eqnarray*}
\chi \left(  \prod_{s=1}^{q} \Delta_{s},\ldots, \prod_{s=1}^{q} \Delta_{q(m-1)+s} \right) &=& \sum_{\sigma \in {NC}_{\geq 2}\left(\left[mq \right], \pi^{*}  \right) }\prod_{V \in \sigma}\mu\left(A_{b_1^{V,\sigma}}\right).
\end{eqnarray*}
Inserting the right hand side of the above equation in (\ref{prespeicher}), we obtain
\begin{eqnarray*}
\kappa_{m}\left( I_{q}^{\hat{N}}(f)\right)  &=& \sum_{i_{1}^{(1)},\ldots,i_{q}^{(1)}}\cdots \sum_{i_{1}^{(m)},\ldots,i_{q}^{(m)}} c_{i_{1}^{(1)}\cdots i_{q}^{(1)}}\cdots c_{i_{1}^{(m)}\cdots i_{q}^{(m)}}\sum_{\sigma \in {NC}_{\geq 2}\left(\left[mq \right], \pi^{*}  \right) }\prod_{V \in \sigma}\mu\left(A_{b_1^{V,\sigma}} \right) \\
&=& \sum_{\sigma \in {NC}_{\geq 2}\left(\left[mq \right], \pi^{*}  \right) }\sum_{i_{1}^{(1)},\ldots,i_{q}^{(1)}}\cdots \sum_{i_{1}^{(m)},\ldots,i_{q}^{(m)}} c_{i_{1}^{(1)}\cdots i_{q}^{(1)}}\cdots c_{i_{1}^{(m)}\cdots i_{q}^{(m)}}\prod_{V \in \sigma}\
\mu\left(A_{b_1^{V,\sigma}} \right).
\end{eqnarray*}
To conclude the proof of \eqref{e:dp}, it remains to notice that 
\begin{equation*}
\sum_{i_{1}^{(1)},\ldots,i_{q}^{(1)}}\cdots \sum_{i_{1}^{(m)},\ldots,i_{q}^{(m)}} c_{i_{1}^{(1)}\cdots i_{q}^{(1)}}\cdots c_{i_{1}^{(m)}\cdots i_{q}^{(m)}}\prod_{V \in \sigma}\mu\left(A_{b_1^{V,\sigma}}\right) = \int_{Z^{\vert \sigma \vert}}f_{\sigma}d\mu^{\vert \sigma \vert}.
\end{equation*}
Relation \eqref{e:dp2} is an immediate consequence of \eqref{e:momcum}. To show \eqref{e:srp}, we exploit 
\eqref{e:dp2} to deduce the crude bound
\[
\varphi\left((I^{\hat{N}}_q(f))^{2m} \right)^{1/2m} \leq  |NC(2mq)|^{1/2m} \times \max\{1 ; KD\}^{q/2}, 
\]
and the desired estimate follows from \eqref{e:stir} and the fact that \[\rho(I^{\hat{N}}_q(f)) = \lim_{m\to\infty} \varphi\left((I^{\hat{N}}_q(f))^{2m} \right)^{1/2m}.\] We have therefore concluded the proof of the Point (i) in the statement for every $f\in \mathscr{E}_q^{00}$. The extension to a general $f\in \mathscr{E}_q$ is achieved by combining the previous computations with Lemma \ref{l:taborn} and Lemma \ref{l:lot}.

\medskip 

\noindent {\it Case 2: ${\mathfrak{M}} = S$.} We start once again with a kernel $f\in \mathscr{E}_q^{00}$. In this case, because the cumulants of the centered semicircular distribution are all zero except for the one of order two, one can rewrite (\ref{sumsurnc}) as 
\begin{equation}
\label{sumdanslecassemicircular}
\kappa \left(  \prod_{s=1}^{q} \Delta_{s},\ldots, \prod_{s=1}^{q} \Delta_{q(m-1)+s} \right) = \sum_{\sigma \in {NC}_{2}\left(\left[mq \right], \pi^{*}  \right) }\prod_{V \in \sigma}\kappa_{2}\left(\Delta_{b_1^{V,\sigma}}\right).
\end{equation}
Note that if the product $mq$ is odd, the above quantity is zero, since in this case ${NC}_{2}\left(\left[mq \right], \pi^{*}  \right) = \emptyset$. Focusing again on (\ref{prespeicher}), equations (\ref{sumsurpartitions})--(\ref{sumdanslecassemicircular}) yield
\begin{equation*}
\kappa_{m}\left( I_{q}^{S}(f)\right)  = \sum_{i_{1}^{(1)},\ldots,i_{q}^{(1)}}\cdots \sum_{i_{1}^{(m)},\ldots,i_{q}^{(m)}} c_{i_{1}^{(1)}\cdots i_{q}^{(1)}}\cdots c_{i_{1}^{(m)}\cdots i_{q}^{(m)}}\sum_{\sigma \in {NC}_{2}\left(\left[mq \right], \pi^{*}  \right) }\prod_{V \in \sigma}\kappa_{2}\left(\Delta_{b_1^{V,\sigma}}\right).
\end{equation*}
Using the fact that, for any $A \in \mathcal{Z}_{\mu}$, $\kappa_{2}\left(S\left(A\right) \right) = \mu\left(A \right) $, we deduce that
\begin{eqnarray*}
\kappa_{m}\left( I_{q}^{S}(f)\right)  &=& \sum_{i_{1}^{(1)},\ldots,i_{q}^{(1)}}\cdots \sum_{i_{1}^{(m)},\ldots,i_{q}^{(m)}} c_{i_{1}^{(1)}\cdots i_{q}^{(1)}}\cdots c_{i_{1}^{(m)}\cdots i_{q}^{(m)}}\sum_{\sigma \in {NC}_{2}\left(\left[mq \right], \pi^{*}  \right) }\prod_{V \in \sigma}\mu\left(A_{b_1^{V,\sigma}} \right) \\
&=& \sum_{\sigma \in {NC}_{2}\left(\left[mq \right], \pi^{*}  \right) }\sum_{i_{1}^{(1)},\ldots,i_{q}^{(1)}}\cdots \sum_{i_{1}^{(m)},\ldots,i_{q}^{(m)}} c_{i_{1}^{(1)}\cdots i_{q}^{(1)}}\cdots c_{i_{1}^{(m)}\cdots i_{q}^{(m)}}\prod_{V \in \sigma}\mu\left(A_{b_1^{V,\sigma}} \right).
\end{eqnarray*}
Using the relation 
\begin{equation*}
\sum_{i_{1}^{(1)},\ldots,i_{q}^{(1)}}\cdots \sum_{i_{1}^{(m)},\ldots,i_{q}^{(m)}} c_{i_{1}^{(1)}\cdots i_{q}^{(1)}}\cdots c_{i_{1}^{(m)}\cdots i_{q}^{(m)}}\prod_{V \in \sigma}\mu\left(A_{b_1^{V,\sigma}} \right) = \int_{Z^{mq/2}}f_{\sigma}d\mu^{mq/2},
\end{equation*}
we deduce \eqref{e:ds}. As before, the moment formula \eqref{e:ds} follows from \eqref{e:momcum}, whereas the spectral bound \eqref{e:srs} is a consequence of the relation
\[
\rho(I^{S}_q(f)) = \lim_{m\to\infty} \varphi\left((I^{S}_q(f))^{2m} \right)^{1/2m},
\]
whose right hand side has to be evaluated according to the sharp arguments provided in \cite[Theorem 5.3.4]{bianespeicher}. The extension to a general $f\in \mathscr{E}_q$ follows from Lemma \ref{l:lot} and from the density of $\mathscr{E}_q^{00}$ in the space $L^2(Z^q)$. 
\qed

\end{document}